\documentclass[11pt]{article}
\usepackage[export]{adjustbox} 
\usepackage{fancyhdr}         
\usepackage{hyperref}        
\usepackage{datetime}    
\usepackage{array}            
\usepackage{amsmath,amssymb}   
\DeclareMathOperator*{\argmin}{arg\,min}
\usepackage[bottom]{footmisc}
\usepackage{mathtools}        
\usepackage{commath}         
\usepackage{enumitem}

\usepackage{caption}         
\usepackage[includeheadfoot,
            left=0.8in,
            right=0.8in,
            top=0.20in,
            bottom=0.20in,
            headheight=8pt]{geometry}
\usepackage{comment}
\usepackage{float}
\usepackage{authblk}

\usepackage{booktabs}
\usepackage{siunitx}
\renewcommand{\headrulewidth}{0pt}
\renewcommand{\footrulewidth}{0pt}

\usepackage[backend=biber,style=numeric,sorting=none,
  giveninits=true,     
  maxbibnames=99,
  date=year,
  doi=false,isbn=false,url=false
]{biblatex}

\DeclareNameAlias{author}{family-given}
\DeclareDelimFormat{multinamedelim}{\addcomma\space}
\DeclareDelimFormat{finalnamedelim}{\addspace\&\space}
\DeclareFieldFormat{pages}{#1}
\newbibmacro*{journalname}{%
  \iffieldundef{shortjournal}
    {\printtext{\mkbibemph{\printfield{journaltitle}}}}
    {\printtext{\mkbibemph{\printfield{shortjournal}}}}}
\DeclareBibliographyDriver{article}{%
  \usebibmacro{author}%
  \setunit{\labelnamepunct}\newblock
  \printfield[title]{title}%
  \newunit\newblock
  \usebibmacro{journalname}%
  \setunit{\space}%
  \printfield[volume]{volume}%
  \setunit{\space}%
  \printtext[parens]{\printfield{year}}%
  \setunit{\addcomma\space}%
  \printfield{pages}%
  \finentry}
\DeclareFieldFormat[article,inproceedings,incollection,periodical]{volume}{\mkbibbold{#1}}
\DeclareFieldFormat{doilink}{\href{https://doi.org/#1}{https://doi.org/#1}}
\DeclareBibliographyDriver{book}{%
  \usebibmacro{author/editor}%
  \setunit{\labelnamepunct}\newblock
  \printtext{\mkbibemph{\printfield{title}}}%
  \iffieldundef{series}{}{%
    \setunit{\addcomma\space}\printfield{series}%
    \iffieldundef{volume}{}{%
      \addcomma\space \printfield{volume}}}%
  \setunit{\addcomma\space}%
  \printtext[parens]{%
    \printlist{publisher}\addcomma\space\printlist{location}\addcomma\space\printfield{year}}%
  \iffieldundef{doi}{}{%
    \par\noindent\printfield[doilink]{doi}}%
  \finentry
}
\DeclareBibliographyDriver{thesis}{%
  \usebibmacro{author}%
  \setunit{\labelnamepunct}\newblock
  \printtext{\mkbibemph{\printfield{title}}}%
  \setunit{\addcomma\space}%
  \printfield{type}
  \setunit{\addcomma\space}%
  \printlist{institution}%
  \setunit{\space}%
  \printtext[parens]{\printfield{year}}%
  \finentry}
\AtEveryBibitem{\clearfield{month}\clearlist{language}}

\begin{document}

\pagestyle{fancy}
\fancyhf{}                   
\fancyfoot[C]{\thepage}        
\fancypagestyle{plain}{%
  \fancyhf{}%
  \fancyfoot[C]{\thepage}%
  \renewcommand{\headrulewidth}{0pt}%
  \renewcommand{\footrulewidth}{0pt}%
}

\title{Data-Driven Modeling of Global Bifurcations and Chaos in a Mechanical System under Delayed and Quantized Control}
\author[1]{Giacomo Abbasciano}
\author[3]{Balázs Endrész}
\author[2]{Gábor Stépán}
\author[1,$\dagger$]{George Haller}
\affil[1]{Department of Mechanical and Process Engineering, ETH Zurich, 8092 Zurich, Switzerland}
\affil[2]{Department of Applied Mechanics, Budapest University of Technology and Economics, H-1111 Budapest, Hungary}
\affil[3]{HUN-REN-BME Dynamics of Machines Research Group, Budapest University of Technology and Economics, H-1111 Budapest, Hungary}

\date{}  
\maketitle
\begingroup
\renewcommand\thefootnote{\textdagger}
\footnotetext{Corresponding author: \href{mailto:georgehaller@ethz.ch}{georgehaller@ethz.ch}}
\endgroup
\let\oldclearpage\clearpage
\let\clearpage\relax

\section*{Abstract}
We illustrate  how the recent theory of Spectral Submanifolds (SSM) can capture global bifurcations and complex dynamics in mechanical systems even under delay and spatial discretization. Specifically, we build a parameter-dependent SSM-reduced model that predicts global heteroclinic and local bifurcations in a Furuta pendulum under control with delay, and verify these predictions numerically. Under additional spatial discretization of the digital controller, we also obtain an SSM-reduced model that correctly reproduces a numerically and experimentally observed microchaotic attractor in the system.
\newpage
\section{Introduction}

Stabilization of unstable equilibria and orbits via feedback control is a critically important task in engineering and applied sciences. Examples include the regulation of robots, brake systems, satellites, and airplanes. Incentive for detailed local stability studies is that systems controlled to operate at an otherwise unstable equilibrium are more agile and 
can execute abrupt departure from an equilibrium with less control effort \cite{Insperger2015, Ott1990}.

Under digital feedback control, however, two challenges arise in stability analysis: feedback delay and spatial discretization. Time delay emerges due to the finite speed of data processing and information transmission, leading to control systems described by delay-differential equations (DDEs) in an infinite-dimensional phase space. Additionally, spatial discretization affects the control action whenever a digital controller is employed, as analog-to-digital converters (ADC) and digital-to-analog converters (DAC) can only handle a finite number of bits in the conversion process.

Beyond these implementational challenges, the physical systems to be controlled often exhibit significant nonlinear behavior. Such nonlinearities can lead to the coexistence of multiple isolated invariant sets such as hyperbolic fixed points, periodic orbits, quasi-periodic tori, and chaotic attractors. Rendering the system globally nonlinearizable, these features cannot be captured simultaneously by linear models \cite{Cenedese2022}, even though control approaches commonly focus on approximating only the linear dynamics \cite{Kaundinya2025}.

For equation-driven settings, powerful tools have been developed to analyze stability of equilibria of systems with time delay. Notably, the semi-discretization method introduced by \cite{Insperger2015,Insperger2002} provides accurate predictions of stability maps for fixed points under PD control \cite{Szaksz2022,Vizi2024}, and has been successfully applied in several engineering applications \cite{InspergerStepan2011}. Nonetheless, the method faces limitations when dealing with high-dimensional systems and parameter uncertainties. As a more recent development, \cite{Szaksz2025} formally computes Spectral Submanifolds (SSMs) in the infinite-dimensional phase space induced by time delays, showing excellent agreement with numerical simulations of the full delayed systems.

Equation-driven approaches, however, become inapplicable when the governing equations are not exactly known, as often occurs, for example, in practice for soft robots and rigid bodies with feedback delay. In such cases, data-driven reduced modeling methodologies become essential.

To address the limitations of linear data-driven model reduction methods \cite{Schmid2010, Williams2015} in capturing non-linearizable dynamics, such as isolated coexisting attractors, recent work has focused on Spectral Submanifold (SSM) reduction \cite{Haller2025}. SSMs are low-dimensional invariant manifolds in the phase space of a dynamical system that are tangent to spectral subspaces identified from linearizing the system at a fixed point. The reduced dynamics of SSMs act as nonlinear reduced order models with which nearby trajectories synchronize exponentially fast.

A rigorous equation-based theory of SSMs is available in several different settings with a broad range of applications across autonomous \cite{Cenedese2022,  Haller2025Videos, Liu2024, Yang2025, Axs2023, Kaszs2024} and nonautonomous systems \cite{Haller2025, Haller2024, Axs2024}. The SSMLearn algorithm developed in \cite{Cenedese2022} enables the identification of SSM-based low-dimensional nonlinear models from data. It has been successfully applied to various physical systems including fluid sloshing experiments \cite{Axs2024}, fluid-structure interaction problems \cite{Xu2024}, pipe flows \cite{Kaszs2024}, non-smooth jointed structures \cite{Morsy2025}, and continuum robots in model predictive control settings \cite{Alora2022, Alora2023, Alora2025}. SSM-reduction has also been extended to piecewise-smooth systems, both autonomous and forced, in \cite{Bettini2024}, and has been rigorously justified for systems governed by delay-differential equations \cite{BuzaHaller2025}. Recent advances, such as the use of Padé approximants in place of polynomial parametrizations \cite{Kaszs2025}, and the development of nonlinear oblique projections onto SSMs \cite{Bettini2025}, have further enhanced the domain of validity of data-driven SSM reconstruction.

An open question has been whether data-driven SSM-reduction can be carried out under parameter dependence. Even more importantly, can a parametric data-driven SSM model predict chaotic dynamics and global bifurcations? Additionally, can all these  be achieved under additional challenges such as the presence of delay and non-smoothness in a physical system? Here we provide a positive answer to these open questions in the context of a specific mechanical system, the Furuta pendulum that is feedback-controlled at its unstable equilibrium \cite{Insperger2015, Vizi2024, Vizi2021, Stepan2017}.

The outline of this paper is as follows. Section \ref{section3} delineates the problem setup for general mechanical systems under delay, when a PD controller with spatial discretization is employed to stabilize a fixed point. 
Section \ref{section6} implements the methodology using data from numerical simulations and experiments on the Furuta pendulum.

\section{Problem setup}\label{section3}
Here we describe the elements we will apply in the SSM-reduced modeling of the controlled Furuta pendulum using data from numerical simulations and experiments.
\subsection{Uncontrolled system}

Consider the autonomous, $n$-dimensional, uncontrolled dynamical system,
\begin{equation}\label{eq:uncontrolledSys}
    \dot x = f(x)=Ax + f_0(x), \quad x\in\mathbb{R}^n, \quad A\in\mathbb{R}^{n\times n}, \quad f_0 = \mathcal{O}(|x|^2)\in C^r, \quad r\in \mathbb{N}^+ \cup \{ \infty, a \} ,
\end{equation}
where $C^a$ refers to the space of analytic functions. Assume that $x=0$ is a hyperbolic fixed point of system~\eqref{eq:uncontrolledSys}.

Let $E \subset \mathbb{R}^n$ denote a spectral subspace of $A$, i.e., a direct sum of a set of real eigenspaces of $A$.
Under appropriate non-resonance conditions on the spectrum of $A$, $E$ is known to have a unique nonlinear continuation, $\mathcal{W}(E)$, in the full nonlinear system~\eqref{eq:uncontrolledSys} (see \cite{Haller2025} for a review of related mathematical results). 
Usually referred to as a (primary) Spectral Submanifold (SSM), $\mathcal{W}(E)$ is a $C^r$ invariant manifold emanating from the origin, has the same dimensions as $E$, and is tangent to $E$ at $x=0$.
The SSM $\mathcal{W}(E)$ is called a slow SSM if all eigenvalues of $A$ outside $E$ have smaller real parts than those inside $E$. 
Furthermore, $\mathcal{W}(E)$ is called a like-mode SSM if the real parts of the eigenvalues inside $E$ have the same sign. Otherwise, $\mathcal{W}(E)$ is called a mixed-mode SSM. The internal dynamics of slow SSMs provide exact reduced-order models with which the dynamics of the full nonlinear system~\eqref{eq:uncontrolledSys} synchronizes exponentially fast \cite{Haller2025,Haller2023}. By the Takens' embedding theorem \cite{Takens1981}, under nondegeneracy conditions on the SSM dynamics, a diffeomorphic copy of a $d$-dimensional SSM $\mathcal{W}(E)$ can be generically embedded in the space of $p>2d$ time-delayed observations of scalar quantity along trajectories of~\eqref{eq:uncontrolledSys}. SSMs can be extracted from observable data using the SSMLearn algorithm developed in \cite{Cenedese2022}, as summarized in \cite{Haller2025}.
We summarize the relevant mathematical results on SSMs and the main steps of the SSMLearn algorithm in Appendix \ref{SectionAppA}.

In this work, system~\eqref{eq:uncontrolledSys} corresponds to the first-order formulation of the equations of motion for the Furuta pendulum~\eqref{eq:eqOfMotion}, which serves as our numerical and experimental test case. The approach we follow, however, applies to data from any mechanical systems seeked to be stabilized around a fixed point using PID control with delayed feedback and spatial discretization, as outlined in Appendix \ref{SectionAppB}.

\subsection{PD controlled system}

We assume that the fixed point at the origin is unstable and seek to stabilize it via feedback control. We start by defining an observable
\[
\xi:\mathbb{R}^n \to \mathbb{R}^l,
\quad
\xi(0)=0,
\quad 
\xi \in C^1,
\]

\noindent
with the help of which we define the tracking error
\[
e(t) = \xi(0) - \xi(x(t)) = -\,\xi(x(t)),
\]
which vanishes at the desired fixed point, $x=0$. We seek to use a PD controller with matrix-valued gains
\[
K_{\rm P}, \; K_{\rm D} \;\in\; \mathbb{R}^{n \times l}.
\]
\noindent
The control input is constructed from the classical proportional and derivative actions:  

\begin{itemize}
\item \textbf{Proportional action:}  $u_P(t) = K_{\rm P}\,e(x(t)) = -\,K_{\rm P}\,\xi(x(t)).$

\item \textbf{Derivative action:} $u_D(t) = K_{\rm D}\,\frac{\rm d}{{\rm d}t}  e(x(t))
     = -\,K_{\rm D}\,\frac{\rm d}{{\rm d}t} \,\xi\bigl(x(t)\bigr).$

\end{itemize}

\noindent

The $n$-dimensional closed-loop dynamics satisfy
\begin{equation}
\dot{x}(t)
=
    f\bigl(x(t)\bigr) + u(x(t)).
\end{equation}
where $u(t) = u_P(t) +  u_D(t) = -\,K_{\rm P}\,\xi(x(t)) - K_{\rm D}\,\frac{\rm d}{{\rm d}t} \,\xi(x(t))
  \;=\; u(x(t))$. The right-hand side depends solely on the state vector $x$ and contains no explicit time dependence. Therefore, the closed-loop system is autonomous, with the origin $x = 0$ as the target equilibrium to be stabilized.

\subsection{PD controlled system with time-delayed feedback}\label{sec:subsec2.3}

Next, we incorporate a time delay $\tau>0$ into the feedback loop. The control input then takes the form
\[
u(t)
= -\,K_{\rm P}\,\xi\bigl(x(t-\tau)\bigr)
  -\,K_{\rm D}\,\frac{\rm d}{{\rm d}t} \,\xi\bigl(x(t-\tau)\bigr).
\]
Accordingly, the closed-loop dynamics are governed by the delay differential equation
\begin{equation}
\dot{x}(t)
=
f\bigl(x(t)\bigr)
- K_{\rm P}\,\xi\bigl(x(t-\tau)\bigr)
- K_{\rm D}\,\frac{\rm d}{{\rm d}t} \,\xi\bigl(x(t-\tau)\bigr),
\label{eq:delayPID}
\end{equation}
which defines a delayed dynamical system with an equilibrium at $x=0$. Unlike finite-dimensional ordinary differential equations (ODEs), delay differential equations (DDEs), such as~\eqref{eq:delayPID}, exhibit dynamics that explicitly depend on the past. Accordingly, their phase space is infinite-dimensional \cite{Diekmann1995} and is given by the the Banach space of continuous history functions
\[
C([-\tau,0],\mathbb{R}^n),
\qquad
\|\varphi\|_\infty = \max_{\theta\in[-\tau,0]} \|\varphi(\theta)\|,
\]
where each state is represented by the history segment
\[
x_t : [-\tau,0]\to \mathbb{R}^n,
\qquad
x_t(\theta) = x(t+\theta), \quad \theta \in [-\tau,0].
\]

To obtain an infinite-dimensional autonomous ODE formulation, we use an alternative representation of~\eqref{eq:delayPID} that can be obtained by introducing past copies of the extended state vector $x$ at integer multiples of the delay. Specifically, we define
\begin{equation}\label{eq:Extension}
x_{\mathrm{ext}}(t)
=
\begin{pmatrix}
x_0(t) \\[2pt]
x_1(t) \\[2pt]
x_2(t) \\[2pt]
\vdots
\end{pmatrix},
\qquad
x_k(t) = x(t-k\tau),
\quad
x_k \in \mathbb{R}^{\,n},
\quad
k \in \mathbb{N},
\end{equation}
Let $e_k(t) = -\xi(x_k(t))$ denote the tracking error at delay layer $k$. The delayed PD law acting on layer $k$ then takes the form
\[
u_k(t)
= K_{\rm P}\,e_{k+1}(t)\;+\;K_{\rm D}\,\frac{\rm d}{{\rm d}t} \,e_{k+1}(t).
\]
Accordingly, the dynamics of the $k^{\rm th}$ layer evolve as
\begin{equation}
\dot{x}_k(t)
=
\underbrace{f\bigl(x_k(t)\bigr)}_{=\,F_k\bigl({x}_k(t)\bigr)}
\;+\;
\underbrace{-K_{\rm P}\,\xi\bigl(x_{k+1}(t)\bigr)
    -K_{\rm D}\,\dot{\xi}\bigl(x_{k+1}(t)\bigr)}_{=\,G_k\bigl(x_{k+1}(t)\bigr)}.
\end{equation}
Thus, in compact notation,
\begin{equation}
\dot{x}_k(t)
= F_k\bigl({x}_k(t)\bigr) + G_k\bigl(x_{k+1}(t)\bigr),
\end{equation}
for $k=0,1,2,\dots$.  
Note that each layer $k$ depends on the subsequent layer $k+1$. Finally, let us introduce the mappings
\[
F(x_{\mathrm{ext}}) \;=\;
\begin{pmatrix}
F_0({x}_0(t))\\[4pt]
F_1({x}_1(t))\\[4pt]
F_2({x}_2(t))\\[2pt]
\vdots
\end{pmatrix},
\qquad
G(x_{\mathrm{ext}}) \;=\;
\begin{pmatrix}
G_0(x_1(t))\\[4pt]
G_1(x_2(t))\\[4pt]
G_2(x_3(t))\\[2pt]
\vdots
\end{pmatrix}.
\]

With these definitions, the delay-controlled dynamics admits the compact representation
\begin{equation}\label{eq:DelayEquation}
\dot{x}_{\mathrm{ext}}(t)
= F\bigl(x_{\mathrm{ext}}(t)\bigr)
+ G\bigl(x_{\mathrm{ext}}(t)\bigr).
\end{equation}
This construction provides an autonomous ODE representation~\eqref{eq:DelayEquation} of the delay system ~\eqref{eq:delayPID}, as its right-hand side depends solely on the extended state $x_{\mathrm{ext}}(t)$ and not explicitly on time, with the origin $x_{\mathrm{ext}} = 0$ representing the fixed point to be stabilized.

\subsection{Temporal sampling and Zero-Order Hold}

We now incorporate a further practical aspect into controlling system~\eqref{eq:uncontrolledSys}: temporal sampling.  
In practice, the delay $\tau$ introduced in Section \ref{sec:subsec2.3} is partially caused by the fact that the control input is computed based on sampled values of the observable vector $\xi(x)$. As shown in \cite{Insperger2015}, such a temporal sampling naturally induces delays in the control loop. Specifically, the control action is evaluated at the beginning of each sampling interval, $t_i = i\Delta t$, where $\Delta t$ is the sampling period, and then it is kept constant until the end of that interval.

 Assuming negligible computation time (absorbed into the sampling-induced delay), the control signal is then actuated and held piecewise constant over the subsequent sampling interval. This implementation is referred to as Zero-Order Hold (ZOH) in control theory. The ZOH mechanism can be extended by introducing an integer parameter $r \geq 0$, representing the number of sampling intervals between measurement and actuation. The case $r=0$ corresponds to the standard ZOH configuration without additional feedback delay, whereas $r=1$ accounts for the standard one-period delay, which will be the setting adopted in our analysis. Details about regimes with further delayed actions can be found in \cite{Insperger2015, InspergerStepan2011}.

The observable $\xi(x(t))$ is sampled at every $\Delta t$, and the corresponding PD control action is actuated and held constant within each sampling interval. If the actuation is shifted at the beginning of the subsequent sampling interval due to measurement or processing delays ($r=1$), one introduces a time-periodic delay function
\begin{equation}\label{eq:rho}
  \rho(t)
  = t + r\,\Delta t - \Delta t\,
    \Bigl\lfloor \tfrac{t}{\Delta t}\Bigr\rfloor,
  \qquad r\in\mathbb{N},
\end{equation}
where $\lfloor \cdot \rfloor$ denotes the integer part. The resulting discrete-time controller, which accounts simultaneously for feedback delay~\eqref{eq:delayPID} and delay purely due to ZOH logic, can therefore be expressed as
\begin{equation}\label{eq:realworlddd}
\dot{x}(t)=
f\bigl(x(t)\bigr)
- K_{\rm P}\,\xi \bigl(x(t-\rho(t))\bigr)
- K_{\rm D}\,\dot \xi\bigl(x(t-\rho(t))\bigr).
\end{equation}
This formulation provides a closed-loop model of the system capturing the combined effects of temporal sampling, Zero-Order Hold, and feedback delay.

The average time delay \cite{Insperger2015, Insperger2002}, denoted by $\overline{\tau}$, is defined as
\begin{equation}\label{eq:average_rho}
  \overline{\tau} = \bar{\rho}
  \;=\;
  \frac{1}{\Delta t}
  \int_{0}^{\Delta t} \rho(t)\,\mathrm{d}t
  \;=\;
  \left(r+\tfrac{1}{2}\right)\Delta t \, .
\end{equation}

Dynamical system~\eqref{eq:realworlddd} is defined in an infinite-dimensional phase space. Rigorous existence, uniqueness and smoothness results for SSM are available for delay systems \cite{BuzaHaller2025}, and a formal SSM reduction carried out directly in the infinite-dimensional phase space of delayed systems has been shown to perform effectively in \cite{Szaksz2025}.
Motivated by these findings, in this work we extend this formalism to incorporate additional non-smooth terms in the SSM reduction of delayed systems, emerging with the introduction of ZOH, as the control signal becomes only piecewise smooth, exhibiting a countable set of discontinuities.
This extension has already demonstrated accurate results for non-smooth systems without delay (see, e.g., \cite{Morsy2025}). The underlying rationale is that a data-driven SSM reduction yields the smoothest reduced-order model consistent with the available data, while retaining the dominant dynamical features of the original system. The accuracy of such a model can then be assessed by predicting trajectories not included in the training dataset.

\subsection{PD controlled system with time-delayed feedback and spatial discretization}

We now account for the additional spatial discretization effect in the control loop, arising from quantization of both measurements and actuation. 

Inclusion of this discretization effect \cite{Stepan2017} into~\eqref{eq:realworlddd} yields

\begin{equation}\label{eq:realWorldSD}
\begin{aligned}
\dot{x}(t)
=
f\bigl(x(t)\bigr) +
 h\left\lfloor  
   - K_{\rm P}\,\Bigl\lfloor \tfrac{1}{h}\,\xi \bigl(x(t-\rho(t))\bigr)\Bigr\rfloor
   - K_{\rm D}\,\Bigl\lfloor \tfrac{1}{h}\,\dot \xi\bigl(x(t-\rho(t))\bigr)\Bigr\rfloor
 \right\rfloor,
\end{aligned}
\end{equation}

\noindent where $h \in \mathbb{R}$ denotes the quantization step. 
Equation~\eqref{eq:realWorldSD} models that the control action is computed and applied as an integer multiple of the quantization step $h$. The inner floor operators $\lfloor \cdot \rfloor$ represent analog-to-digital conversion (ADC), while the outer floor operators represent digital-to-analog conversion (DAC). For simplicity, we assume identical quantization resolution $h$ for both ADC and DAC.

Analogously to the construction in~\eqref{eq:Extension}, we extend the discretized system by introducing infinitely many past copies of the state vector. 
We define the $k^{\rm th}$ layer as
\[
x_{k}(t) \\[2pt]
=
x\!\bigl(t - \operatorname{sgn}(k)\,\rho(t)
           - H(k - 2)\,r (k-1)\,\Delta t \bigr) 
\in \mathbb{R}^{\,n},
\quad k=0,1,2,\dots
\]
\noindent where we have used the Heaviside function $H(x)=0$ for $x<0$, $H(x)=1$ for $x\ge0$ and the sign function $\operatorname{sgn}(x)=-1$ for $x<0$, $0$ for $x=0$, $1$ for $x>0$. We also introduce the extended state as
\begin{equation}\label{eq:Extension1}
x_{\mathrm{ext}}(t)
=\begin{pmatrix}
x_{0}(t) \\[2pt]
x_{1}(t) \\[2pt]
x_{2}(t) \\[2pt]
x_{3}(t) \\[2pt]
\vdots
\end{pmatrix}
=
\begin{pmatrix}
x(t) \\[2pt]
x(t-\rho (t)) \\[2pt]
x(t-\rho (t)-r\Delta t) \\[2pt]
x(t-\rho (t)-2r\Delta t) \\[2pt]
\vdots
\end{pmatrix},
\end{equation}
\noindent At delay layer $k$, the dynamics takes the form
\begin{equation}\label{eq:delayed_layers}
\begin{aligned}
\dot x_{k}(t)
=
\underbrace{
f\bigl(x_{k}(t)\bigr)}_{\,F_{{\rm d},k}\bigl(x_{k}(t)\bigr)}
\;+\;
\underbrace{
h \Bigl\lfloor -K_{\rm P}\,\Bigl\lfloor \tfrac{1}{h}\,\xi\bigl(x_{k+1}(t)\bigr)\Bigr\rfloor
- K_{\rm D}\,\Bigl\lfloor \tfrac{1}{h}\,\dot{\xi}\bigl(x_{k+1}(t)\bigr)\Bigr\rfloor \Bigr\rfloor}_{\,G_{{\rm d},k}\bigl(x_{k+1}(t)\bigr)\;=\;\text{input from next delay layer}}.
\end{aligned}
\end{equation}

\noindent
In compact notation,
\begin{equation}
\dot{x}_{k}(t)
= F_{{\rm d},k}\bigl(x_{k}(t)\bigr) + G_{{\rm d},k}\bigl(x_{k+1}(t)\bigr),
\qquad k=0,1,2,\dots,
\end{equation}
showing that each layer depends on the subsequent one. As before, this necessitates an infinite extension of the phase space. For the entire extended system we introduce the mappings
\[
F_{\rm d}(x_{\mathrm{ext}})
=
\begin{pmatrix}
F_{{\rm d},0}(x_{0}(t))\\[4pt]
F_{{\rm d},1}(x_{1}(t))\\[4pt]
F_{{\rm d},2}(x_{2}(t))\\[2pt]
\vdots
\end{pmatrix},
\qquad
G_{\rm d}(x_{\mathrm{ext}})
=
\begin{pmatrix}
G_{{\rm d},0}(x_{1(t)})\\[4pt]
G_{{\rm d},1}(x_{2(t)})\\[4pt]
G_{{\rm d},2}(x_{3(t)})\\[2pt]
\vdots
\end{pmatrix}.
\]
\noindent
We obtain the infinite-dimensional ODE version of ~\eqref{eq:realWorldSD} in the form:
\begin{equation}\label{eq:ODESD11}
\dot{x}_{\mathrm{ext}}(t)
= F_{\rm d}\bigl(x_{\mathrm{ext}}(t)\bigr)
+ G_{\rm d}\bigl(x_{\mathrm{ext}}(t)\bigr),
\end{equation}

\noindent where equation \eqref{eq:ODESD11} models temporal sampling, Zero-Order Hold (ZOH) logic, delay in the feedback, and spatial discretization. Dynamical system \eqref{eq:ODESD11} is autonomous in its infinite-dimensional phase space. One can thus proceed by seeking a diffeomorphic copy of the SSM in a delay embedding space by the classical Takens' theorem. 
The spatial quantization can generate a non-smooth chaotic attractor in the neighborhood of the fixed point of interest, as has been shown in previous studies~\cite{Stepan2017, Haller1996}.  
\section{Application to a feedback-delayed PD-stabilized mechanical system: the Furuta pendulum}\label{section6} 

\begin{figure}[H]
  \centering  \includegraphics[width=0.45\textwidth]{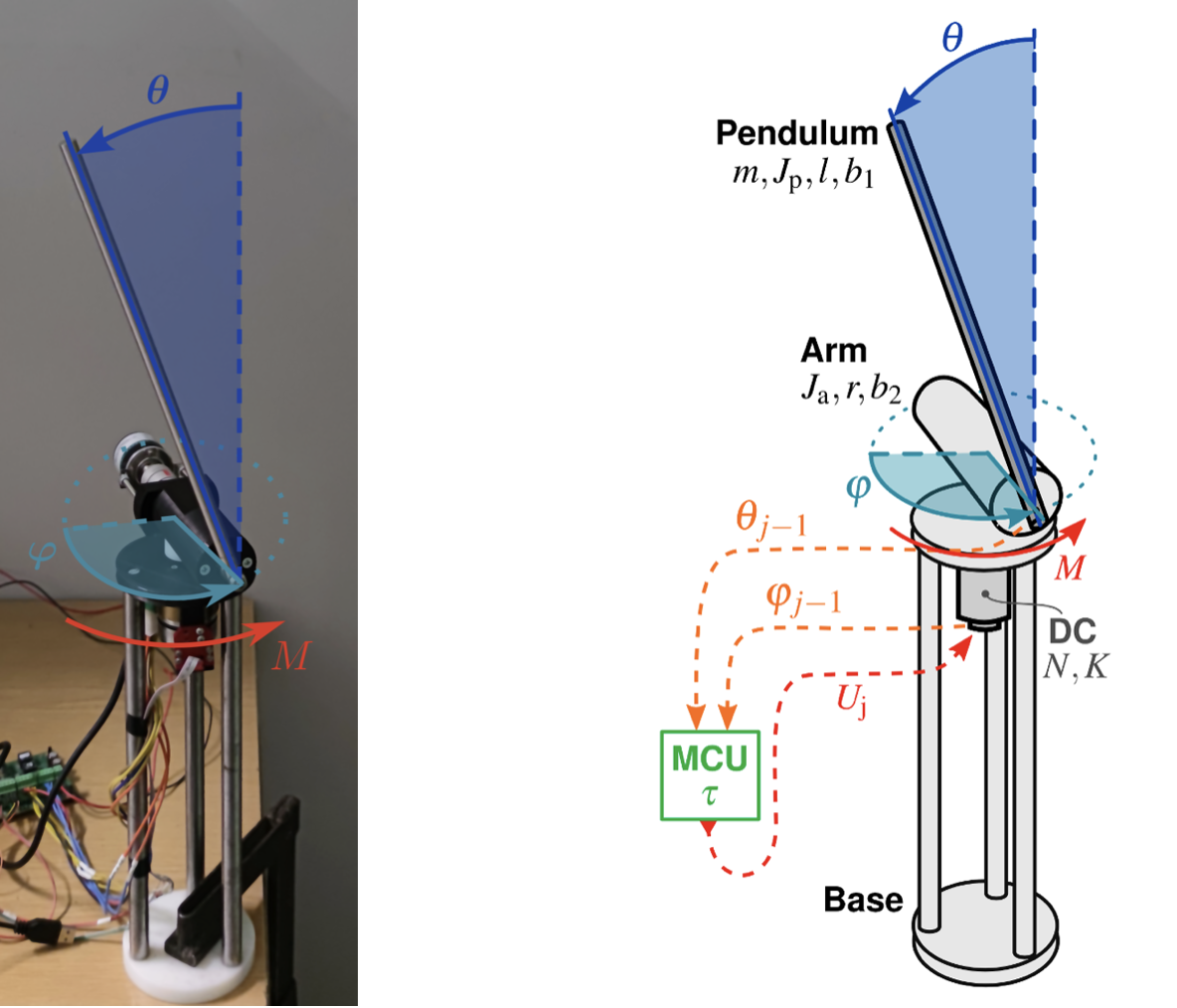}
    \caption{Left: experimental rig used for controlling the Furuta pendulum. Right: schematic of the experimental setup (adopted from \cite{Vizi2024}).}   
  \label{fig:expSetups}
\end{figure}

The Furuta pendulum, also known as the rotary inverted pendulum, is a two-degree-of-freedom (2 DOF) system with angular coordinates $\theta$ and $\varphi$ (see Fig.~\ref{fig:expSetups}). Unlike the planar inverted pendulum, its pivot undergoes circular motion \cite{Vizi2024}. Despite its relative experimental simplicity, analyzing this system poses two key challenges: the interaction between time delay and the non-smoothness introduced by spatial discretization in the stabilization of the nonlinear system. Importantly, the analysis carried out here will apply to any mechanical system whose equilibria are stabilized by PID controllers under time delay and spatial discretization.  

 The right hand sides of the equations of motion of the uncontrolled system are of class $C^\infty$, ensuring the existence of mixed-mode SSMs emanating from the saddle-type equilibrium corresponding to the upright position of the pendulum. For the controlled system, the right-hand side only remains $C^\infty$ within each sampling interval due to the piecewise-smooth control action induced by ZOH logic.

\subsection{Numerical simulations with time delay and no spatial discretization} 

In closed-loop control of nonlinear systems, feedback delays can be critical for stabilization. Beyond a certain delay threshold, the critical time delay \cite{Insperger2015}, the equilibrium cannot be stabilized for any choice of proportional and derivative gains in the two-dimensional $(K_{\rm P}, K_{\rm D})$ parameter space. In this case, no stable region exists in the stability map.  

Stability maps are formed by two bounding curves in the $(K_{\rm P}, K_{\rm D})$ plane: the static loss-of-stability curve and the dynamic loss-of-stability curve. Crossing the former causes a real eigenvalue of the linearized system at the equilibrium to move into the right half of the complex plane. Crossing the latter causes a complex conjugate pair of eigenvalues to cross the imaginary axis and destabilize the system. For sufficiently large delays, the dynamic loss-of-stability curve may self-intersect. At this self-intersection, a codimension-two bifurcation occurs, where two distinct complex conjugate eigenvalue pairs are simultaneously purely imaginary \cite{Vizi2024}. The fixed point then bifurcates into a two-dimensional invariant torus. Depending on the ratio of the imaginary parts of the eigenvalues, the torus may be resonant or non-resonant: in the former case, trajectories close on periodic orbits, whereas in the latter they densely wind on the torus.  

\begin{figure}[H]
    \centering
    \includegraphics[width=0.60\linewidth]{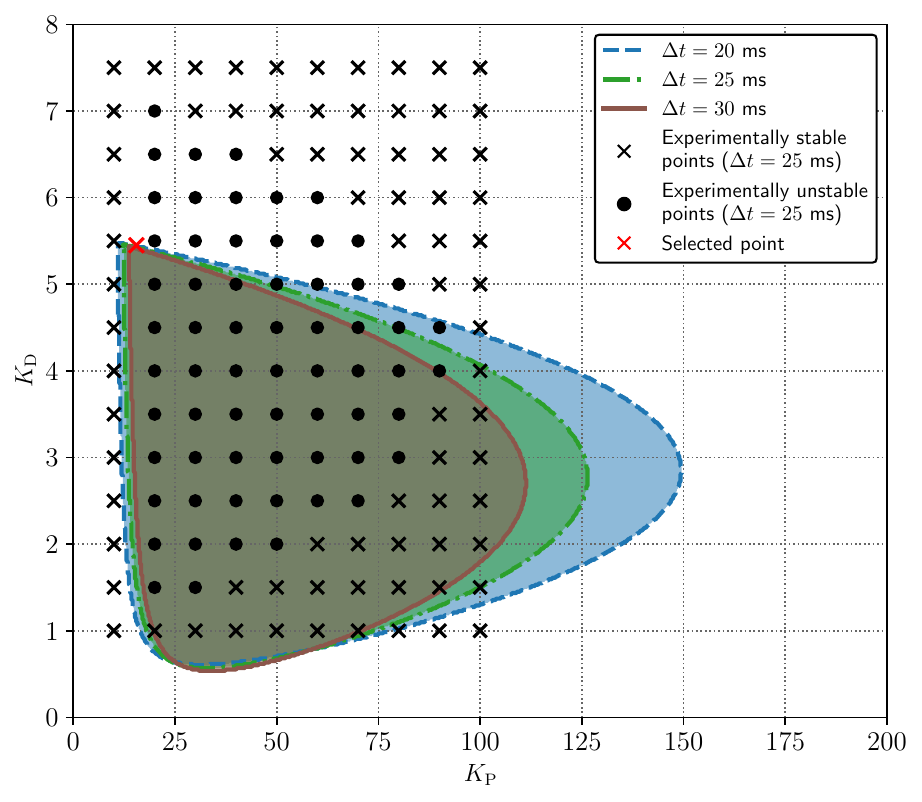}
    \caption{Stability map of the Furuta pendulum in the plane of the control gains. The theoretical stability boundaries are presented for three different sampling delays with the parameters presented in Table~\ref{tab:furutaParams2}. The experimentally observed stable and unstable points are adapted from \cite{Vizi2024}. For numerical simulations, the control gains at the red cross are used, which is slightly outside the stable region, close to the double-Hopf bifurcation point.}
    \label{fig:stabMap}
\end{figure}

For the Furuta pendulum shown in Fig.~\ref{fig:expSetups}, we use the stability map in Figure ~\ref{fig:stabMap} from \cite{Vizi2021} to select the parameter point $(K_{\rm P}, K_{\rm D}) = (15.5\,\mathrm{V}, 5.45\,\mathrm{Vs})$ at $\Delta t=25\,\mathrm{ms}$, located near the self-intersection of the dynamic loss-of-stability curve in the stability map. 
At this point, the upright equilibrium is unstable but lies close to the codimension-two bifurcation point.  In the numerical simulations, we integrate the full nonlinear equations of motion of the Furuta pendulum under ZOH delayed PD feedback $u(t)$, with finite-difference derivatives inducing an effective one-sample delay $(r=1)$:
\begin{equation}
\begin{alignedat}{2}
\dot{\theta} \;&=\;&& \omega_\theta, \\
\dot{\omega}_\theta \;&=\;&& \frac{1}{\Delta(\theta)}\Big[\left(J_{\rm a} + J_{\rm p} \sin^2\theta\right)\left(J_{\rm p} \sin\theta \cos\theta\, \dot{\varphi}^2 - b_1 \omega_\theta + m g l \sin\theta \right) +\\&&&+ m r l \cos\theta\, \left(\mathcal{M}(t)-2J_{\rm p}\sin\theta\cos\theta\dot{\theta}\dot{\varphi}-b_2\dot{\varphi}+mrl\dot{\theta}^2\right) \Big], \\
\dot{\varphi} \;&=\;&& \omega_\varphi, \\
\dot{\omega}_\varphi \;&=\;&& \frac{1}{\Delta(\theta)}\Big[J_{\rm p}\left(\mathcal{M}(t)-2J_{\rm p}\sin\theta\cos\theta\dot{\theta}\dot{\varphi}-b_2\dot{\varphi}+mrl\dot{\theta}^2\right) +\\&&&+ m r l \cos\theta \left(J_{\rm p} \sin\theta \cos\theta\, \dot{\varphi}^2 - b_2 \omega_\theta + m g l \sin\theta \right) \Big],
\end{alignedat}\label{eq:eqOfMotion}
\end{equation}

\noindent where
\begin{equation}
\begin{aligned}
\Delta(\theta) &=
J_{\rm p}\!\left(J_{\rm a} + J_{\rm p}\sin^2\theta\right)
- m^2 r^2 l^2 \cos^2\theta,
\qquad
\mathcal{M}(t) = N\,u(t) - K\,\dot{\varphi}(t), \\
u(t) &=
- K_{\rm P}\,\theta\!\left(t-\rho(t)\right)
- K_{\rm D}\,\frac{\theta\!\left(t-\rho(t)\right)-\theta\!\left(t-\rho(t)-\Delta t\right)}{\Delta t}
+ K_{\varphi {\rm D}}\,\frac{\varphi\!\left(t-\rho(t)\right)-\varphi\!\left(t-\rho(t)-\Delta t\right)}{\Delta t}.
\end{aligned}
\end{equation}
All physical and control parameters are listed in Table~\ref{tab:furutaParams2}.


\newcommand{\pentry}[2]{%
\begin{tabular}[c]{@{}c@{}}
$#1$\\[-1pt]
$#2$
\end{tabular}}

\begin{table}[H]
\centering
\renewcommand{\arraystretch}{1.05}
\begin{tabular}{@{}c@{\hspace{6pt}}c@{\hspace{6pt}}c@{\hspace{6pt}}c@{\hspace{6pt}}c@{\hspace{6pt}}c@{\hspace{6pt}}c@{\hspace{6pt}}c@{\hspace{6pt}}c@{\hspace{6pt}}c@{\hspace{6pt}}c@{\hspace{6pt}}c@{\hspace{6pt}}c@{}}
\toprule
$m$ & $l$ & $r$ & $g$ & $J_{\rm p}$ & $J_{\rm a}$ & $b_1$ & $b_2$ & $N$ & $K$ & $K_{\rm P}$ & $K_{\rm D}$ & $K_{\varphi {\rm D}}$ \\[-2pt]
{\scriptsize kg} &
{\scriptsize m} &
{\scriptsize m} &
{\scriptsize $\mathrm{m\,s^{-2}}$} &
{\scriptsize $\mathrm{kg\,m^{2}}$} &
{\scriptsize $\mathrm{kg\,m^{2}}$} &
{\scriptsize N\,m\,s} &
{\scriptsize N\,m\,s} &
{\scriptsize --} &
{\scriptsize N\,m\,s} &
{\scriptsize V} &
{\scriptsize V\,s} &
{\scriptsize V\,s} \\[-2pt]
0.191 &
0.15 &
0.094 &
9.81 &
$5.73\times 10^{-3}$ &
$1.34\times 10^{-3}$ &
0.039 &
0.02094 &
1.05 &
1.12706 &
15.5 &
5.45 &
1.5 \\
\bottomrule
\end{tabular}
\caption{Physical and control parameters used in the numerical simulations \cite{Vizi2021}.}
\label{tab:furutaParams2}
\end{table}

At this stage of modeling, the spatial discretization was not included. Although digital implementations inherently possess finite bit resolution, here we assume infinite precision in order to study the effect of the time delay on its own.  

From the simulated trajectories, we compute a 4D mixed-mode autonomous SSM in a delay-embedding space, using only time series of angular displacement of the vertical arm $\theta$. The SSM is anchored to the equilibrium corresponding to the upright position (a saddle-type fixed point in the uncontrolled system, mapped to the origin after a coordinate shift).

The minimal possible dimension of an SSM containing the 2-torus arising from the codimension two bifurcation of this system is $d=4$. Such a bifurcation arises when two pairs of complex-conjugate eigenvalues of the linearized system cross the imaginary axis. To construct a parametric reduced-order model, we compute autonomous SSMs together with their normal-form reduced dynamics in polar coordinates for four distinct values of the time delay (equivalently, the sampling time $\Delta t$): $30.5\,\mathrm{ms}$, $31\,\mathrm{ms}$, $32\,\mathrm{ms}$, and $32.5\,\mathrm{ms}$ while mantaining constant $(K_{\rm P}, K_{\rm D}) = (15.5\,\mathrm{V}, 5.45\,\mathrm{Vs})$. These delay values are selected because, starting from $25\,\mathrm{ms}$, a progressive increase in the delay leads to higher degree of instability and the emergence of increasingly rich dynamical behavior.

Here we employ linear interpolation to construct the parametric SSM-reduced model because the full system does not depend smoothly on the delay parameter. Consequently, smooth interpolation methods (e.g., spline, polynomials) are not justified and, in fact, perform worse near bifurcations, where they may introduce artificial oscillations that are not present in the manifold and in the reduced dynamics when in the phase space extended with parameters. Nevertheless, when the underlying system is smooth with respect to the parameter of interest, spline interpolation has been observed to outperform linear interpolation, as expected.

Figure~\ref{fig:phasePortraits1} shows phase portraits of the parametric SSM model in the polar coordinates $(\rho_{1}, \rho_{2})$ plane as $\Delta t$ varies. These results illustrate the dependence of the reduced dynamics on the delay parameter, including predictions at intermediate parameter values not used in the training of the parametric model. 

\newcommand{\imgwidth}{0.5\textwidth}

\begin{figure}[H]
  \centering

  \begin{minipage}[b]{\imgwidth}
    \centering
    \includegraphics[width=\linewidth]{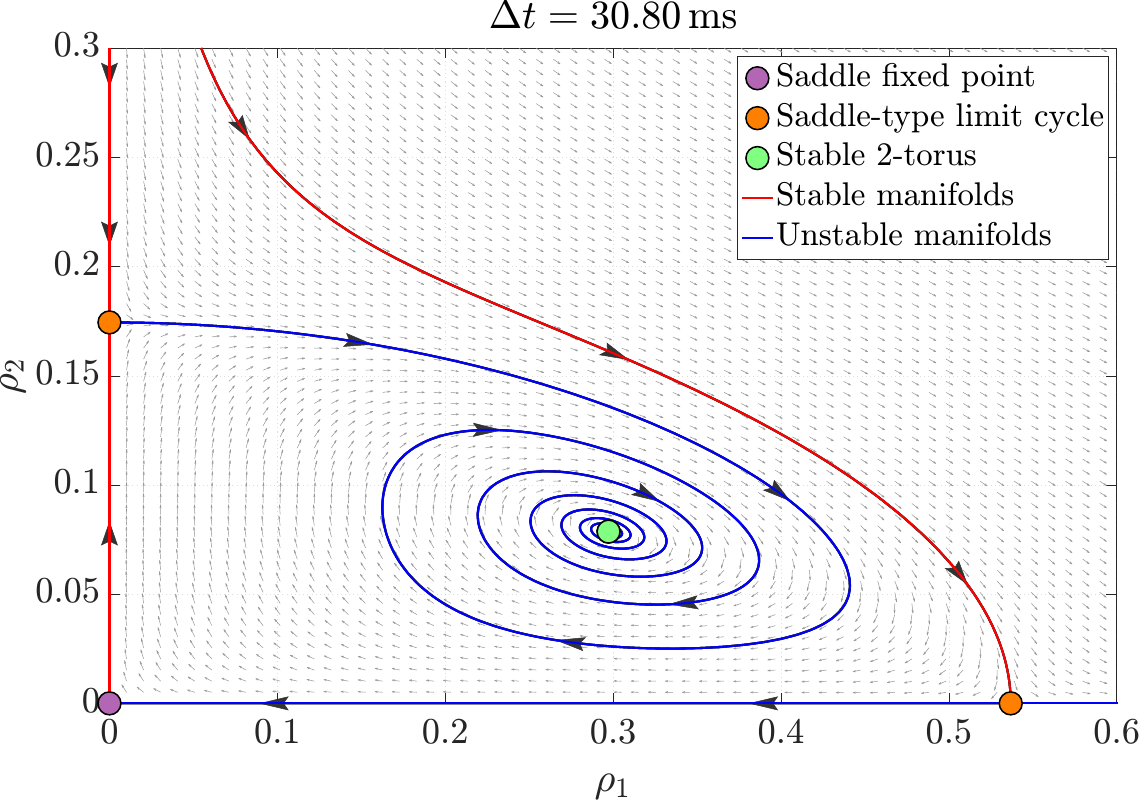}
    \caption*{(a) Sampling time 30.80 ms (32.47 Hz)}
  \end{minipage}\hfill
  \begin{minipage}[b]{\imgwidth}
    \centering
    \includegraphics[width=\linewidth]{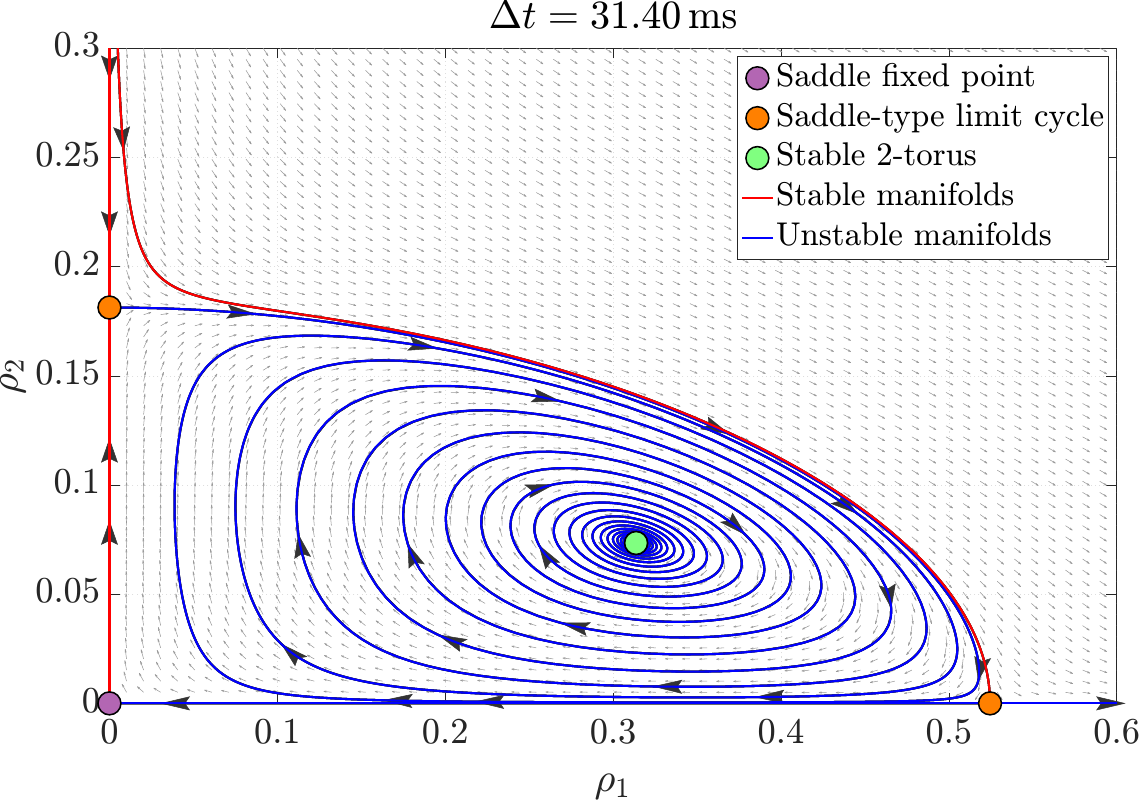}
    \caption*{(b) Sampling time 31.40 ms (31.85 Hz)}
  \end{minipage}
\end{figure}
\begin{figure}[H]
  \centering
  \begin{minipage}[b]{\imgwidth}
    \centering
    \includegraphics[width=\linewidth]{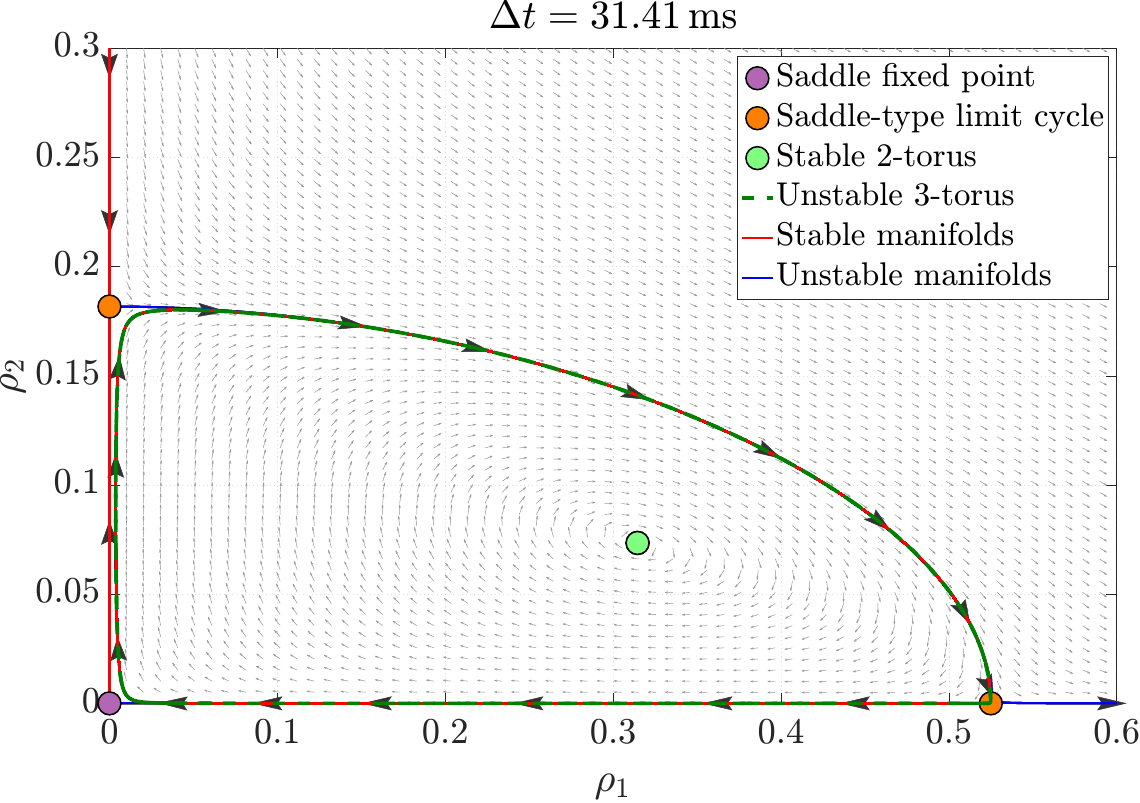}
    \caption*{(c) Sampling time 31.41 ms (31.84 Hz)}
  \end{minipage}\hfill
  \begin{minipage}[b]{\imgwidth}
    \centering
    \includegraphics[width=\linewidth]{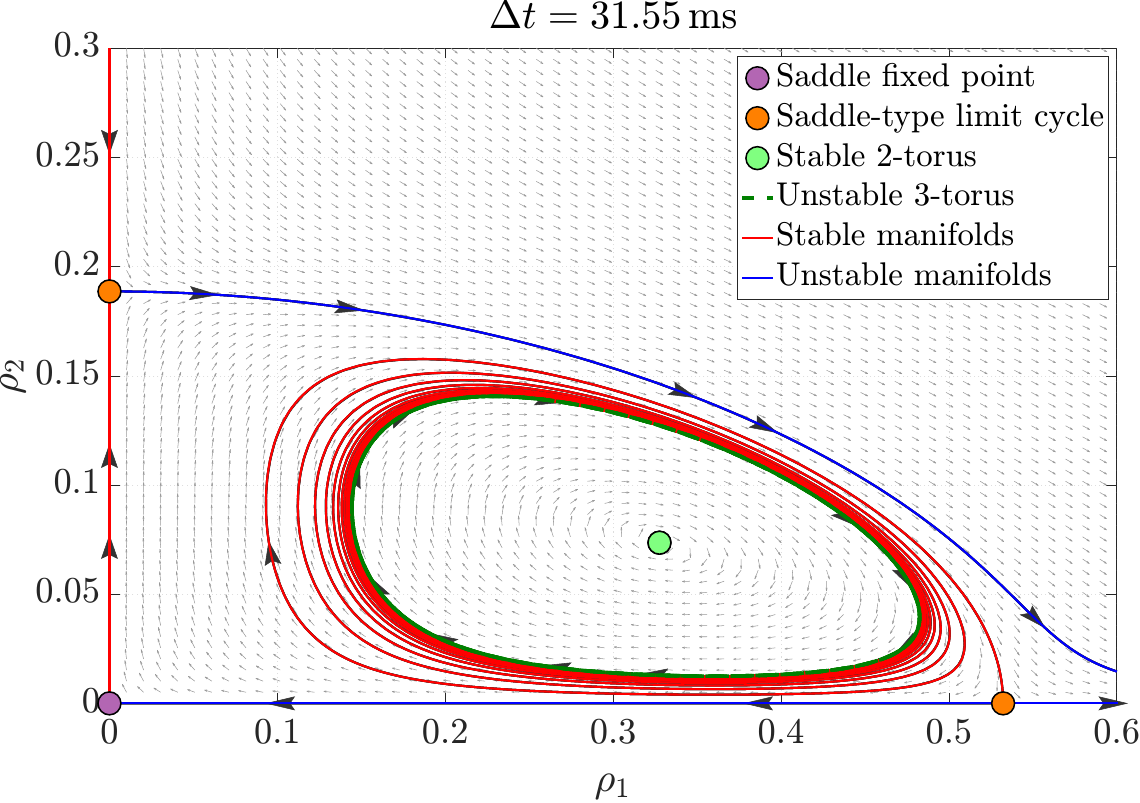}
    \caption*{(d) Sampling time 31.55 ms (31.70 Hz)}
  \end{minipage}
\end{figure}

\clearpage

\begin{figure}[H]
  \centering

  \begin{minipage}[b]{\imgwidth}
    \centering
    \includegraphics[width=\linewidth]{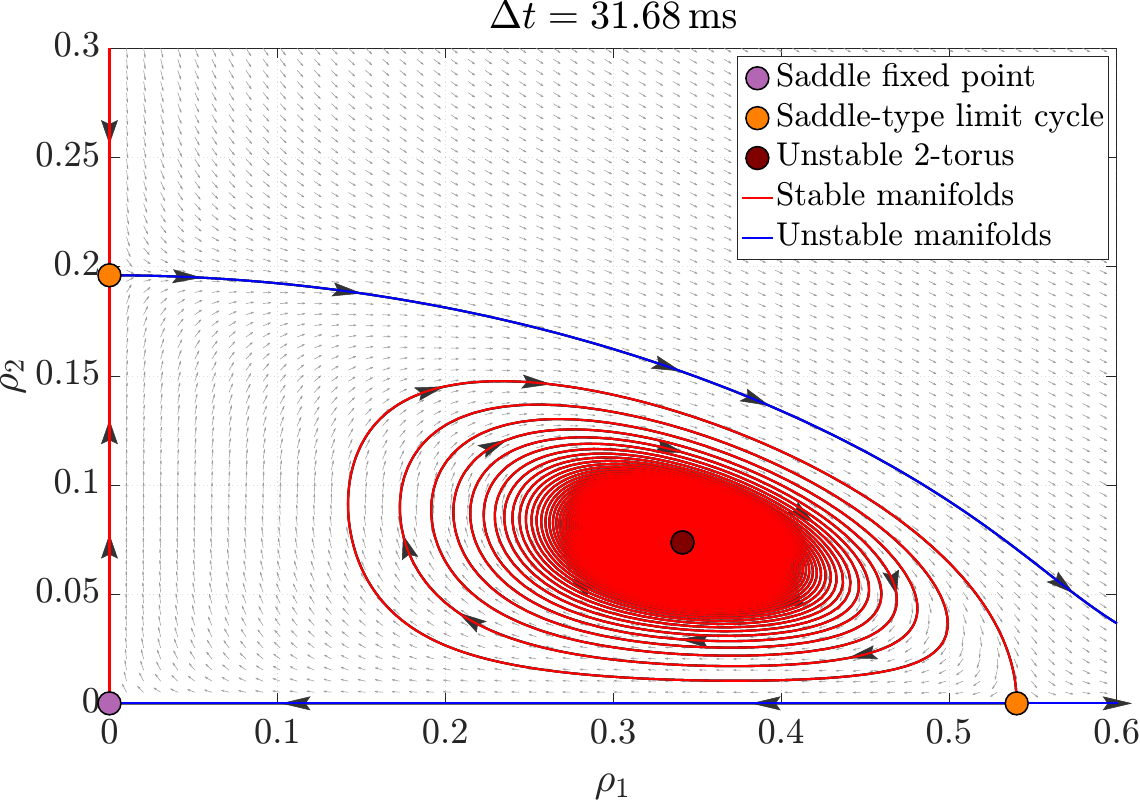}
    \caption*{(e) Sampling time 31.68 ms (31.57 Hz)}
  \end{minipage}\hfill
  \begin{minipage}[b]{\imgwidth}
    \centering
    \includegraphics[width=\linewidth]{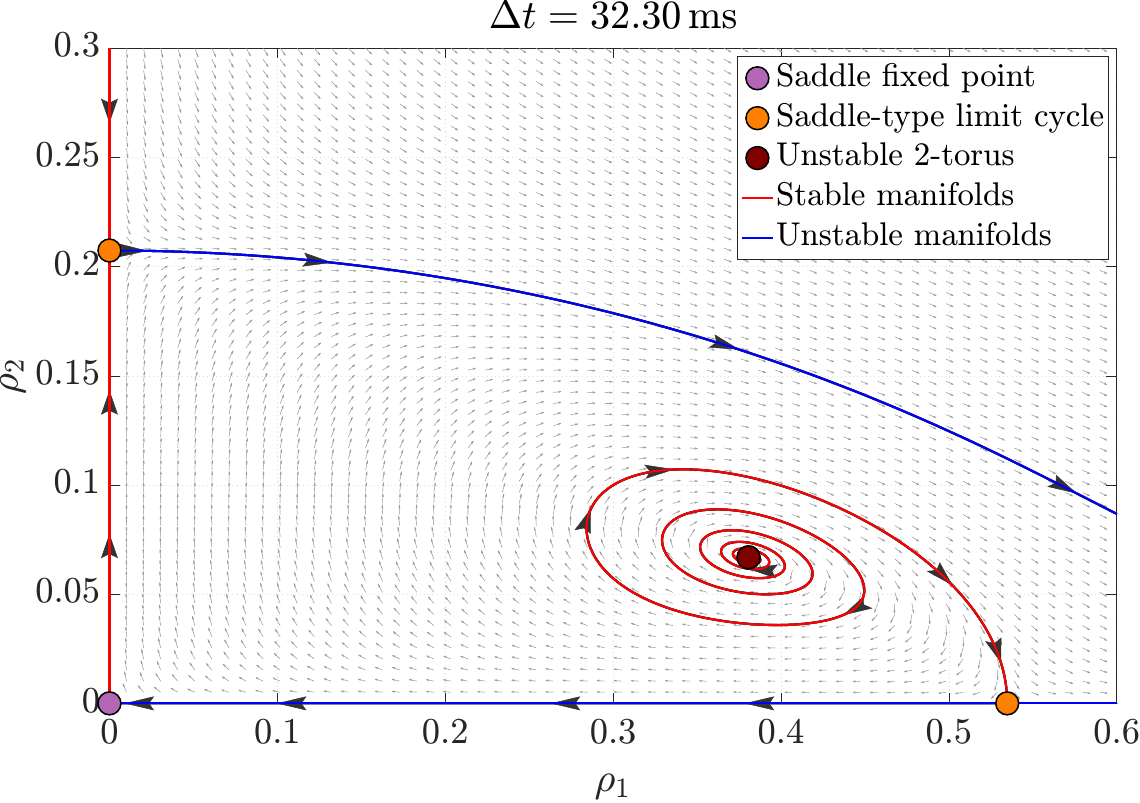}
    \caption*{(f) Sampling time 32.30 ms (30.96 Hz)}
  \end{minipage}

  \caption{Predicted phase portraits of the data-driven parametric SSM‐reduced model at six unseen values of the feedback delay parameter.}
  \label{fig:phasePortraits1}
\end{figure}

The parameter sweep in Fig. ~\ref{fig:phasePortraits1}(a) initially demonstrates the presence of a stable 2-torus, confirming the theoretical predictions from \cite{StepanHaller1995, Molnar2017}, and revealing the existence of two saddle-type limit cycles. The stable and unstable manifolds of these two limit cycles merge into a heteroclinic connection as the time delay in the feedback is increased, as seen in Fig.~\ref{fig:phasePortraits1}(b–c). This connection subsequently breaks under further increases in time delay, resulting in the emergence of an unstable 3-torus through the heteroclinic (global) bifurcation of limit cycles \cite{GuckenheimerHolmes1983}.
The 3-torus provides an escape route for the state toward attracting $\mathcal{O}(1)$ fixed points introduced by the controller, corresponding to pendulum rotation about the horizontal arm with a tilted vertical arm. Increasing the time delay further results in the shrinkage of the 3-torus, as seen in Fig.~\ref{fig:phasePortraits1}(d). This process continues until the stable 2-torus undergoes a subcritical Hopf bifurcation, leading to a change in its stability type. As a consequence, the 2-torus becomes unstable, amidst the disappearance of the 3-torus, as shown in Figs.~\ref{fig:phasePortraits1}(e–f).  

Table~\ref{tab:normal-form-coeffs} reports the coefficients of the truncated normal form governing the reduced amplitude dynamics $\dot{\rho}_1$ and $\dot{\rho}_2$ for the same values of the sampling time $\Delta t$ in Figure~\ref{fig:phasePortraits1}.

\begin{table}[H]
\centering
\footnotesize
\setlength{\tabcolsep}{2pt}
\renewcommand{\arraystretch}{0.9}
\begin{tabular}{c
S S S S S S
@{\hspace{2pt}}
c
S S S S S S}
\toprule
\multicolumn{7}{c}{$\dot{\rho}_1(\rho_1,\rho_2)$} &
\multicolumn{7}{c}{$\dot{\rho}_2(\rho_1,\rho_2)$} \\
\cmidrule(lr){1-7}\cmidrule(lr){8-14}
Monomial &
30.80 & 31.40 & 31.41 & 31.55 & 31.68 & 32.30 &
Monomial &
30.80 & 31.40 & 31.41 & 31.55 & 31.68 & 32.30 \\
\midrule
$\rho_1$          & -0.11 & -0.10 & -0.10 & -0.09 & -0.09 & -0.07 &
$\rho_2$          &  0.27 &  0.27 &  0.27 &  0.27 &  0.26 &  0.27 \\
$\rho_1^3$        &  0.04 &  0.06 &  0.06 &  0.07 &  0.07 &  0.10 &
$\rho_1^2\rho_2$  & -2.12 & -1.99 & -1.99 & -1.87 & -1.77 & -1.65 \\
$\rho_1\rho_2^2$  & 11.65 & 10.76 & 10.71 & 10.03 &  9.41 &  8.21 &
$\rho_2^3$        & -8.21 & -7.43 & -7.39 & -6.90 & -6.44 & -5.88 \\
$\rho_1^5$        &  1.21 &  1.05 &  1.04 &  0.91 &  0.79 &  0.58 &
$\rho_1^4\rho_2$  & -2.16 & -1.62 & -1.60 & -1.25 & -0.93 & -0.09 \\
$\rho_1^3\rho_2^2$& 41.12 & 35.61 & 35.28 & 30.62 & 26.28 & 17.57 &
$\rho_1^2\rho_2^3$&-27.41 &-23.89 &-23.66 &-20.36 &-17.31 & -7.09 \\
$\rho_1\rho_2^4$  &109.98 & 93.75 & 92.78 & 79.24 & 66.66 & 41.66 &
$\rho_2^5$        &-22.49 &-19.71 &-19.43 &-15.52 &-11.88 &-10.85 \\
\bottomrule
\end{tabular}
\caption{Normal form coefficients of $\dot{\rho}_1$ (left) and $\dot{\rho}_2$ (right) for values of the sampling time $\Delta t\,\,\,{\rm [ms]}$ in Fig.~\ref{fig:phasePortraits1}.}
\label{tab:normal-form-coeffs}
\end{table}

\begin{figure}[H]
  \centering
  \captionsetup{justification=centering}

  \begin{minipage}[t]{1\textwidth}
    \centering
    \includegraphics[width=10cm]{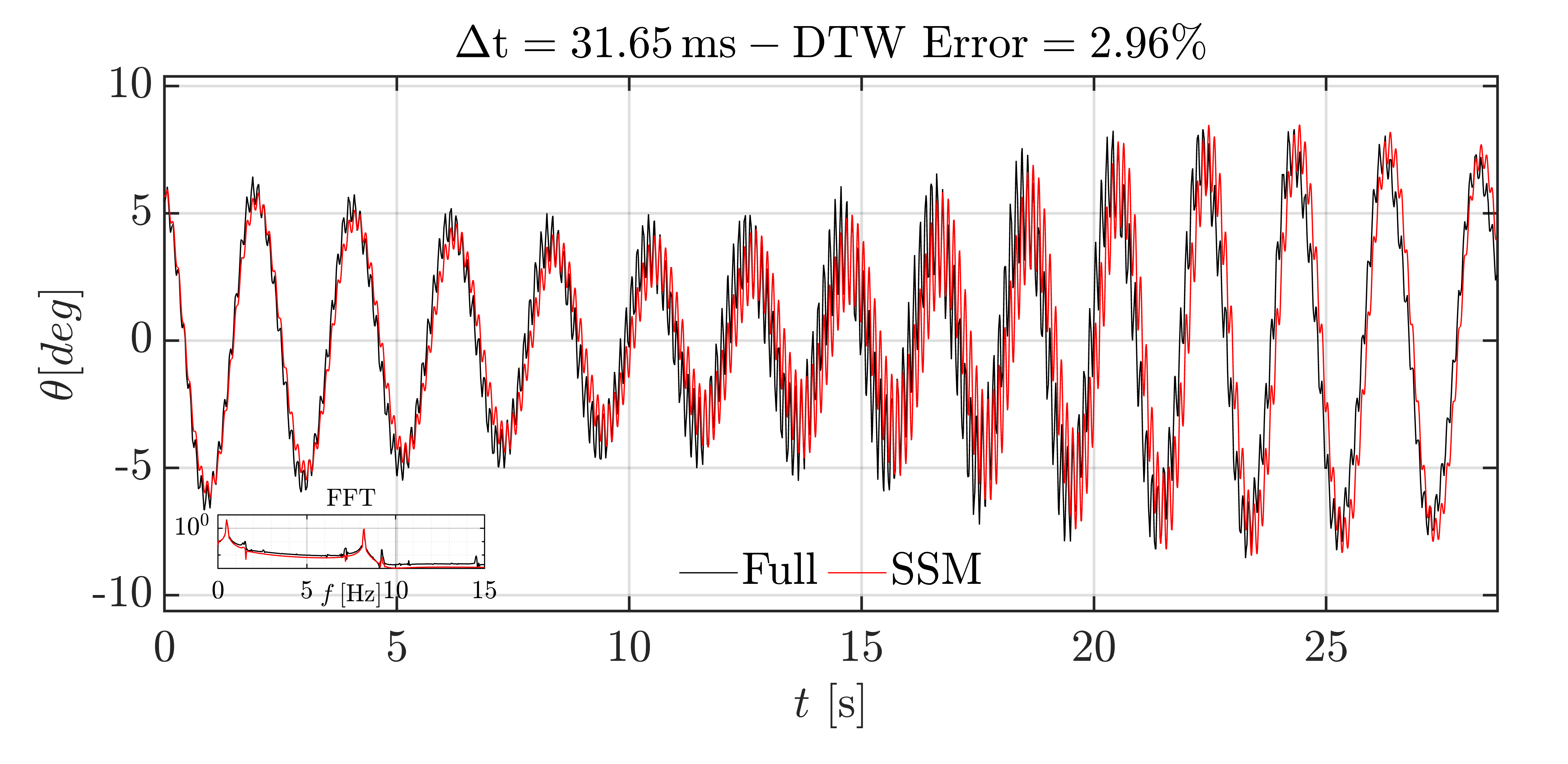}
    \vspace{-8pt}
    \caption*{(f) Sampling time: 31.65\,ms (31.60 Hz)}
  \end{minipage}
  \vspace{-8pt}
  \captionsetup{justification=raggedright,singlelinecheck=false}

  \caption{Predicted trajectories at parameter values not used in the training: comparison between the full system trajectories (black) and the SSM-model trajectories (red). The corresponding FFT spectra are shown in the lower-left inset.}
  \label{fig:TrajectoryPredictions1}
\end{figure}

In Fig.~\ref{fig:TrajectoryPredictions1}, we verify model predictions on test trajectories simulated at parameter values not used in training the data-driven SSM-reduced model.

\begin{figure}[H]
  \centering

  \begin{minipage}[b]{\imgwidth}
    \centering
    \includegraphics[width=\linewidth]{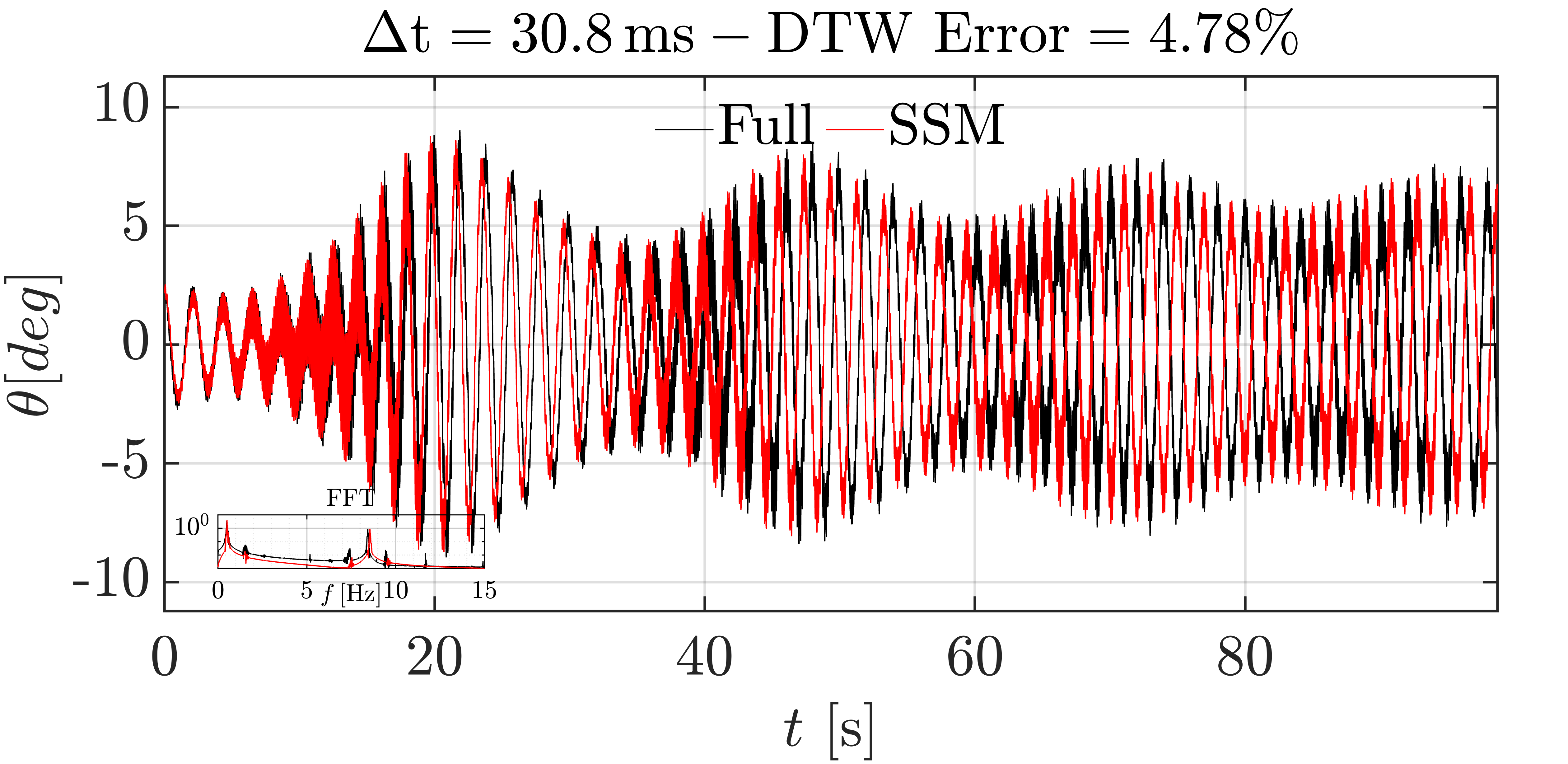}
    \vspace{-10pt}
    \caption*{(a) Sampling time 30.80\,ms (32.47\,Hz)}
  \end{minipage}\hfill
  \begin{minipage}[b]{\imgwidth}
    \centering
    \includegraphics[width=\linewidth]{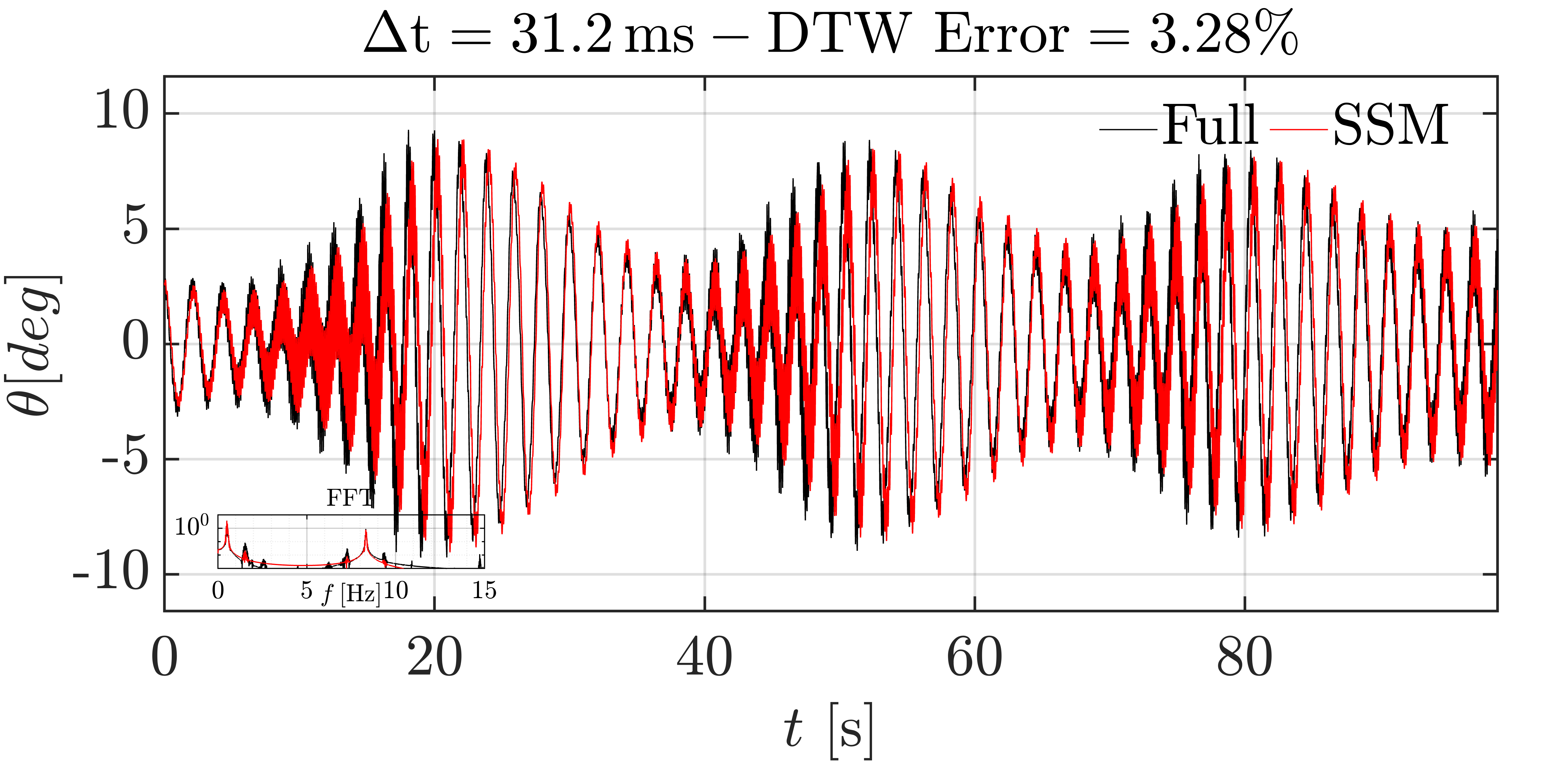}
    \vspace{-10pt}
    \caption*{(b) Sampling time 31.20\,ms (32.05\,Hz)}
  \end{minipage}

\end{figure}
  
\vspace{-10pt}

\begin{figure}[H]
  \centering

  \begin{minipage}[b]{\imgwidth}
    \centering
    \includegraphics[width=\linewidth]{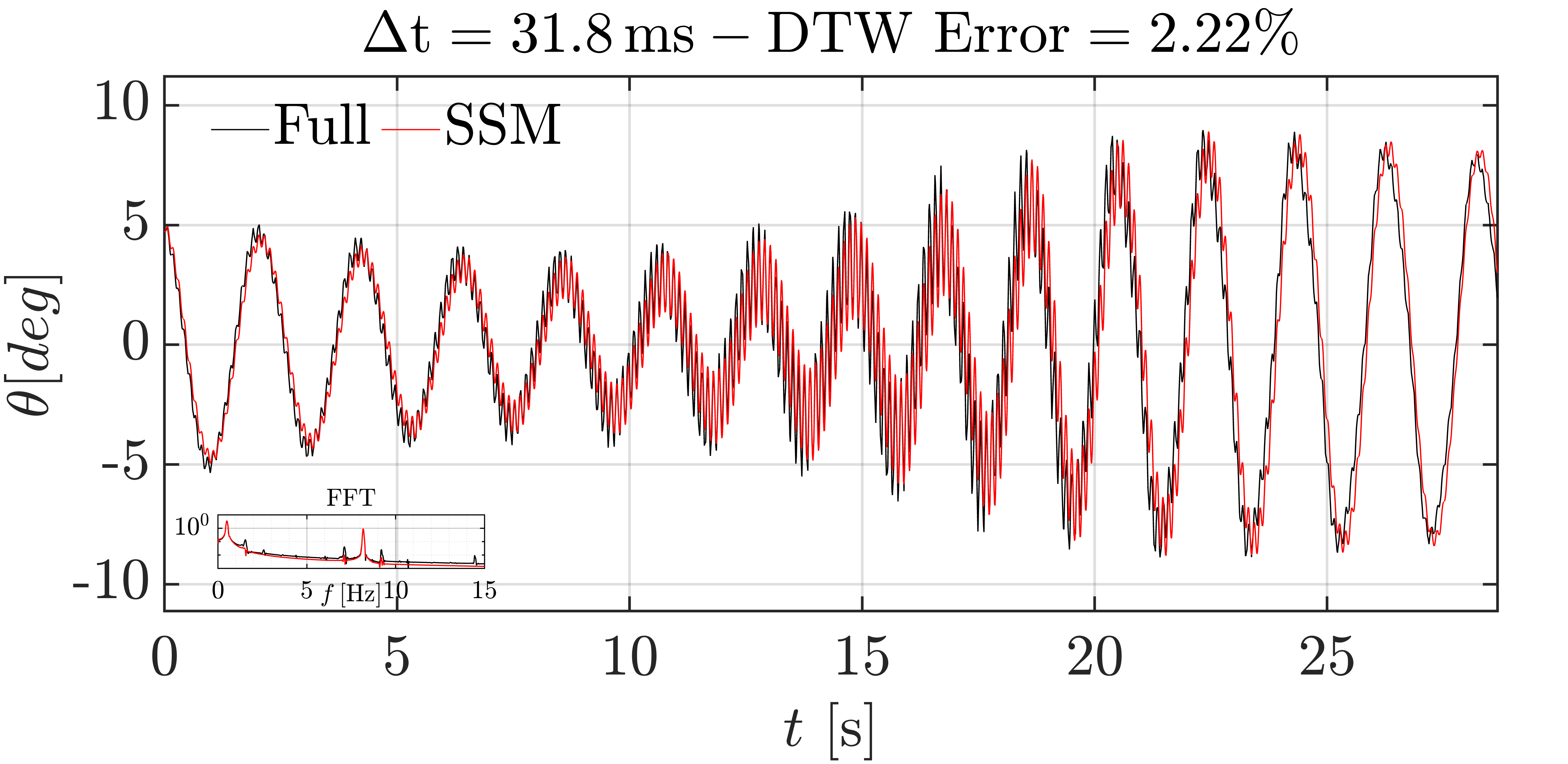}
    \vspace{-10pt}
    \caption*{(c) Sampling time 31.80\,ms (31.45\,Hz)}
  \end{minipage}\hfill
  \begin{minipage}[b]{\imgwidth}
    \centering
    \includegraphics[width=\linewidth]{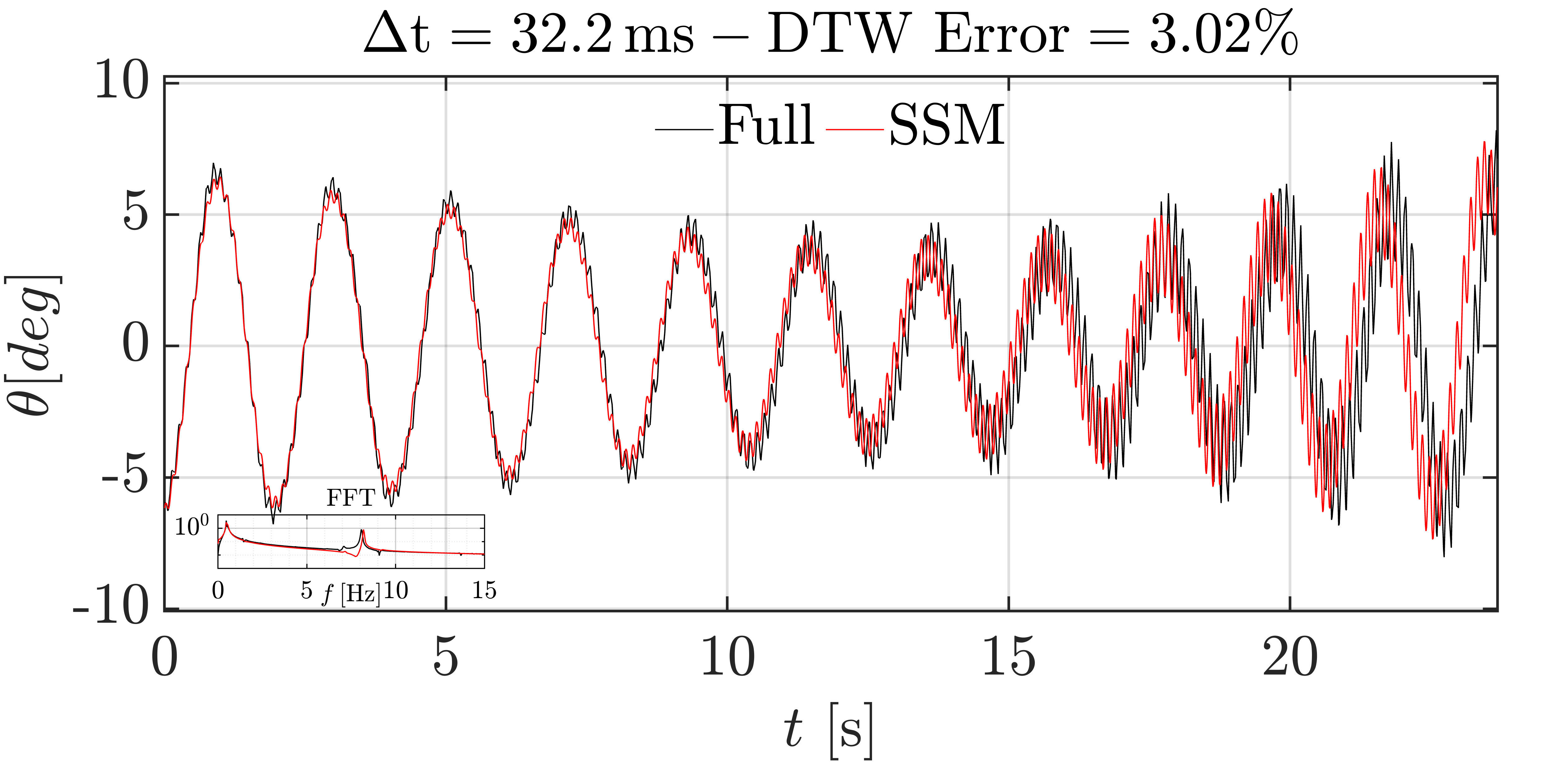}
    \vspace{-10pt}
    \caption*{(d) Sampling time 32.20\,ms (31.06\,Hz)}
  \end{minipage}

\end{figure}

\vspace{-10pt}

\begin{figure}[H]

  \centering
  \captionsetup{justification=centering}

  \begin{minipage}[t]{1\textwidth}
    \centering
    \includegraphics[width=10cm]{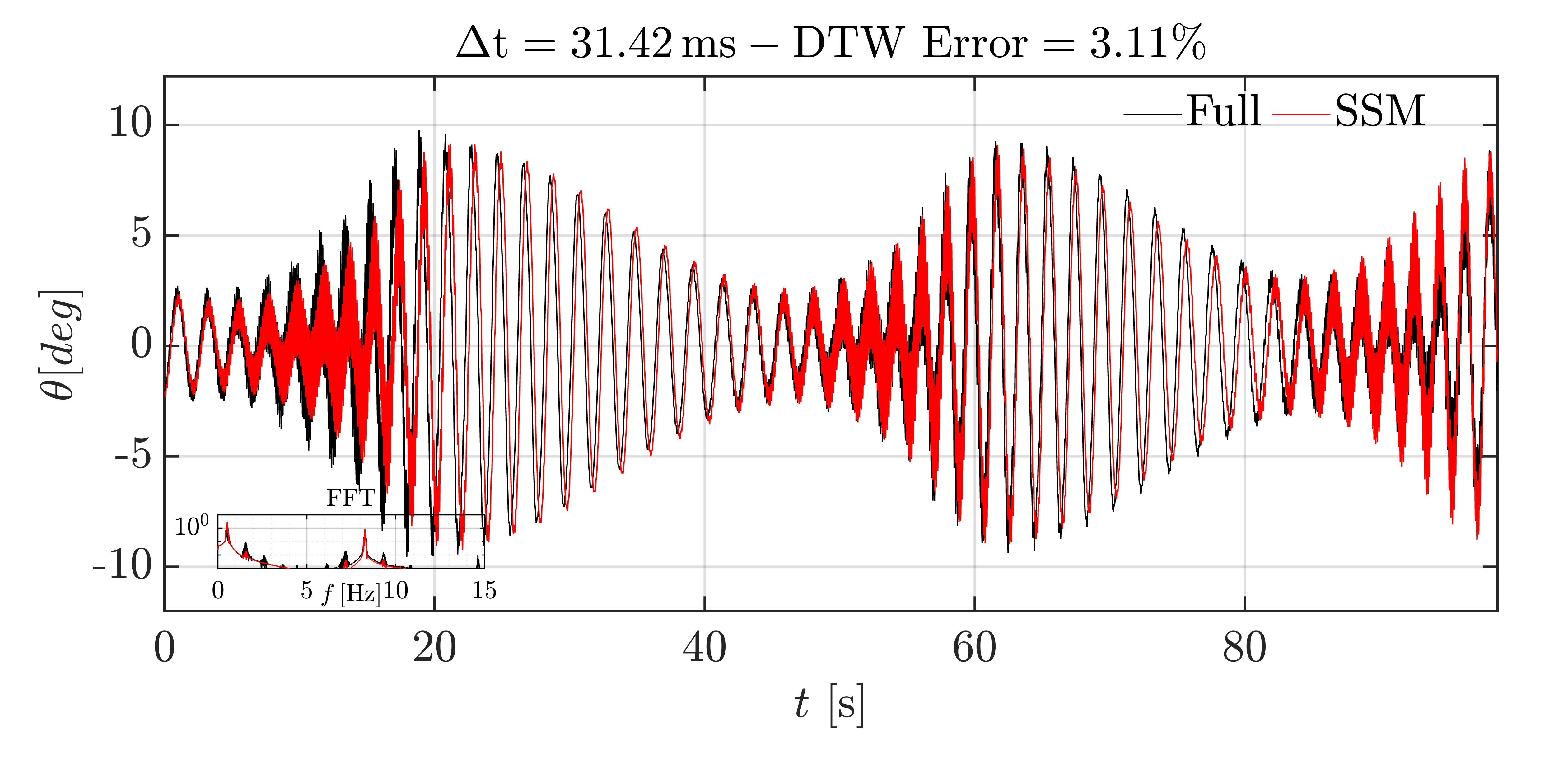}
    \vspace{-8pt}
    \caption*{(e) Sampling time: 31.42\,ms (31.83 Hz)}
  \end{minipage}
\end{figure}

\vspace{-8pt}
  
Phase errors originates from shear in the flow map, which is always present but intensifies near bifurcations, as shown in the phase portraits Fig.~\ref{fig:TrajectoryPredictions1}(e-f). Even small errors in the initial conditions, when projected orthogonally to the tangent space, can amplify into significant phase discrepancies: trajectories complete several revolutions before converging to the stable 2-torus as the parameter approaches the heteroclinic bifurcation.  

A similar phenomenon arises at the Hopf bifurcation: although the 2-torus is unstable, trajectories still complete several revolutions before converging to a fixed point created further away from the origin by the PD controller. Even small modeling inaccuracies, as can be seen by the spectrum plots in ~\ref{fig:TrajectoryPredictions1}, inevitably introduce a time-dependent phase error (at leading order, perturbations to the trajectories $\xi_0 \doteq y_0-x_0$ are advected by the gradient of the flow‐map). For this reason, as long as the phase error is moderate and the model trajectory closely follows the full system's trajectory, the reduced model is considered accurate. In general, one cannot expect more, given the sensitivity to initial conditions typical even in non-chaotic systems.  

To mitigate the dominance of the phase error in our measurements of the model prediction error, we introduce an alternative error measure: the Normalized Mean Trajectory Error (NMTE) defined in \cite{Cenedese2022}, applied to trajectories aligned using Dynamic Time Warping (DTW) \cite{SarinEtAl2010}. This approach isolates the amplitude error while reducing the influence of phase discrepancies.

The domain of attraction of the stable 2-torus predicted by the model is validated by advecting the full system from three distinct initial conditions. The trajectory shown in brown in Fig.~\ref{fig:LALALALA}, initialized within the domain of attraction bounded by the 3-torus, asymptotically converges to the stable 2-torus. In contrast, the two trajectories shown in yellow and blue, initialized outside the domain of attraction, converge to $\mathcal{O}(1)$-distant attracting equilibria. These trajectories pass through the narrow channel formed by the stable and unstable manifolds of the two saddle-type limit cycles, in agreement with the gyroscopic low-frequency loss-of-stability mechanism described in \cite{Haller1992}.

\newcommand{\imgWWwidth}{0.54\textwidth} 
\newcommand{\imgWWWwidth}{0.44\textwidth} 

\begin{figure}[H] 
  \centering
  \begin{minipage}[c]{\imgWWwidth}
    \centering
    \includegraphics[width=\linewidth]{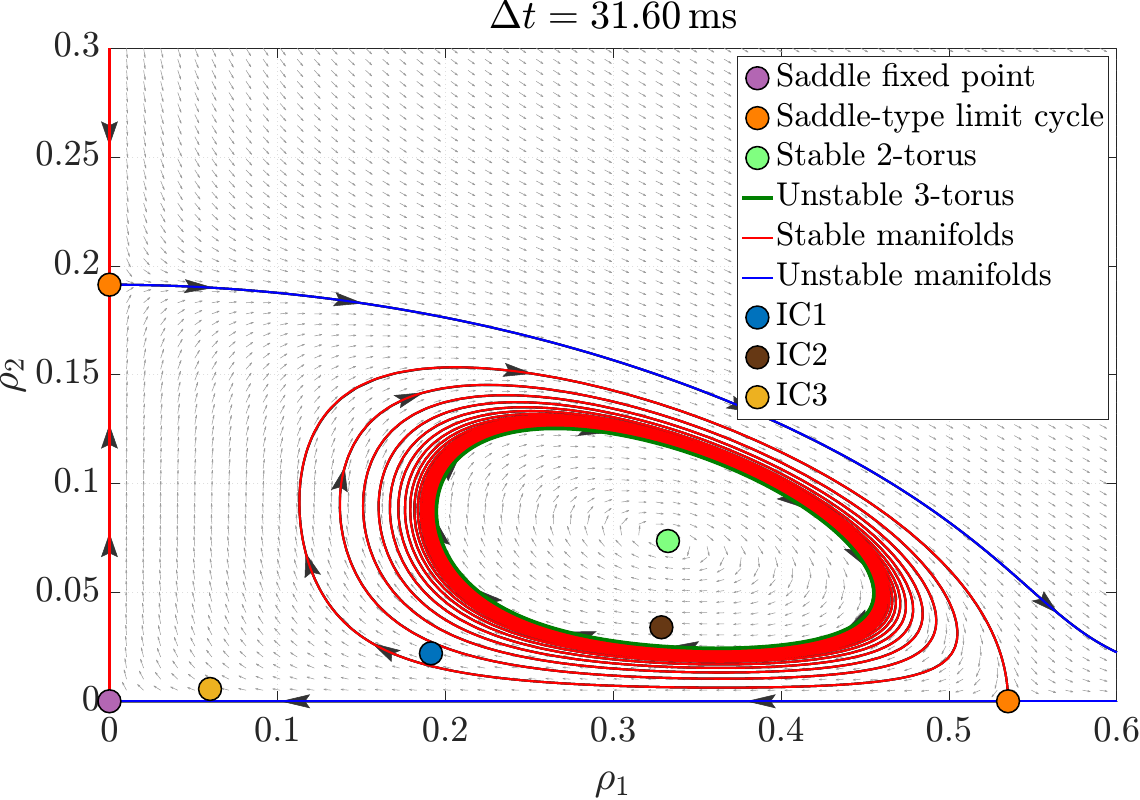}
    \label{fig:mp1}
  \end{minipage}\hfill
  \begin{minipage}[c]{\imgWWWwidth}
    \centering
    \includegraphics[width=\linewidth]{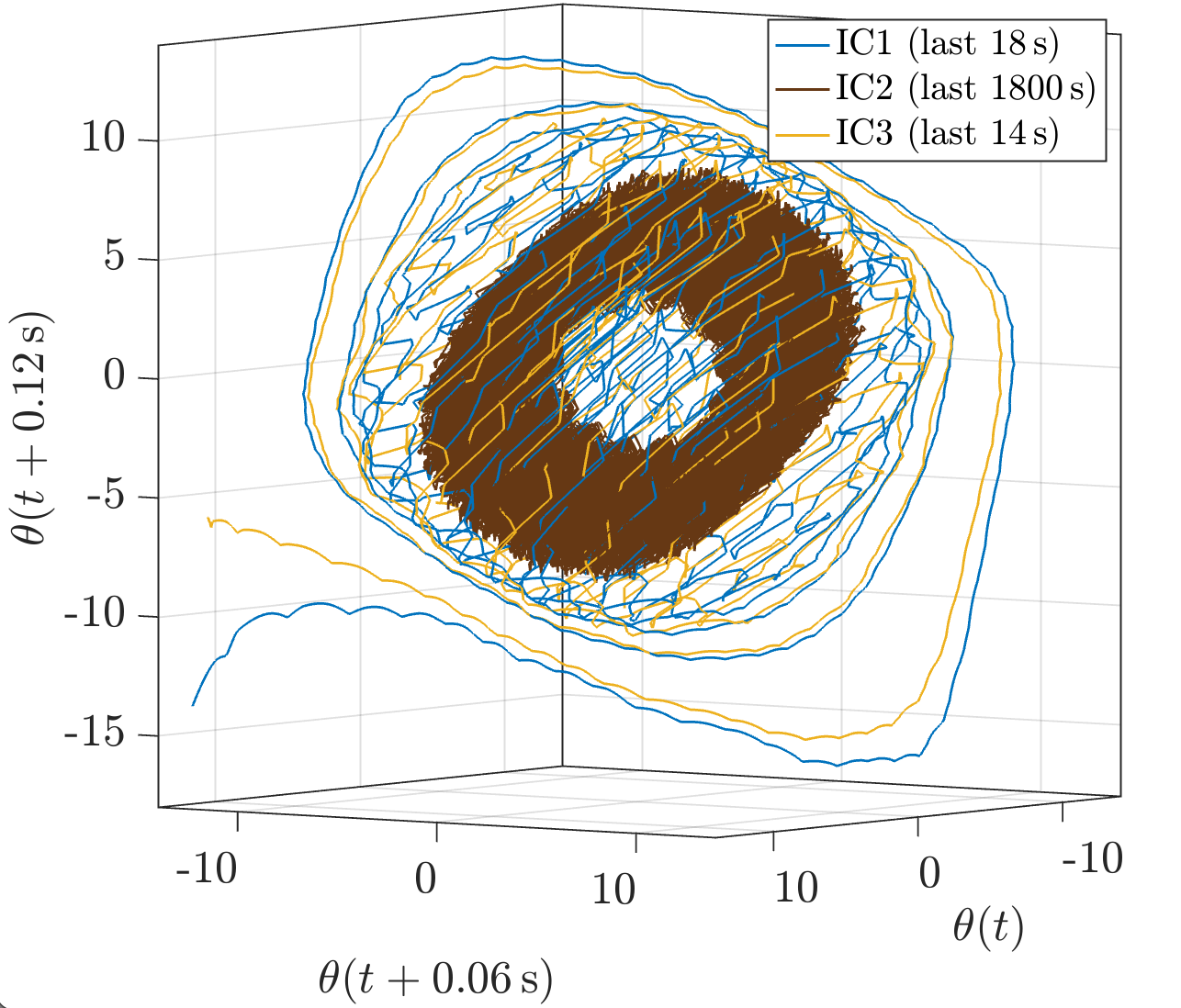}
    \label{fig:mp2}
  \end{minipage}
  \caption{Left: Phase portrait of the parametric SSM model of the Furuta pendulum at a sampling time of $31.60\,\mathrm{ms}$. Three full-system trajectories are initialized at the colored markers (IC1: blue, IC2: brown, IC3: yellow). Right: the same three full-system trajectories shown along time and in a 3D projection of the delay-embedding space.}
  \label{fig:LALALALA}
\end{figure}

\begin{figure}[H]
  \centering
  \includegraphics[width=0.55\textwidth]{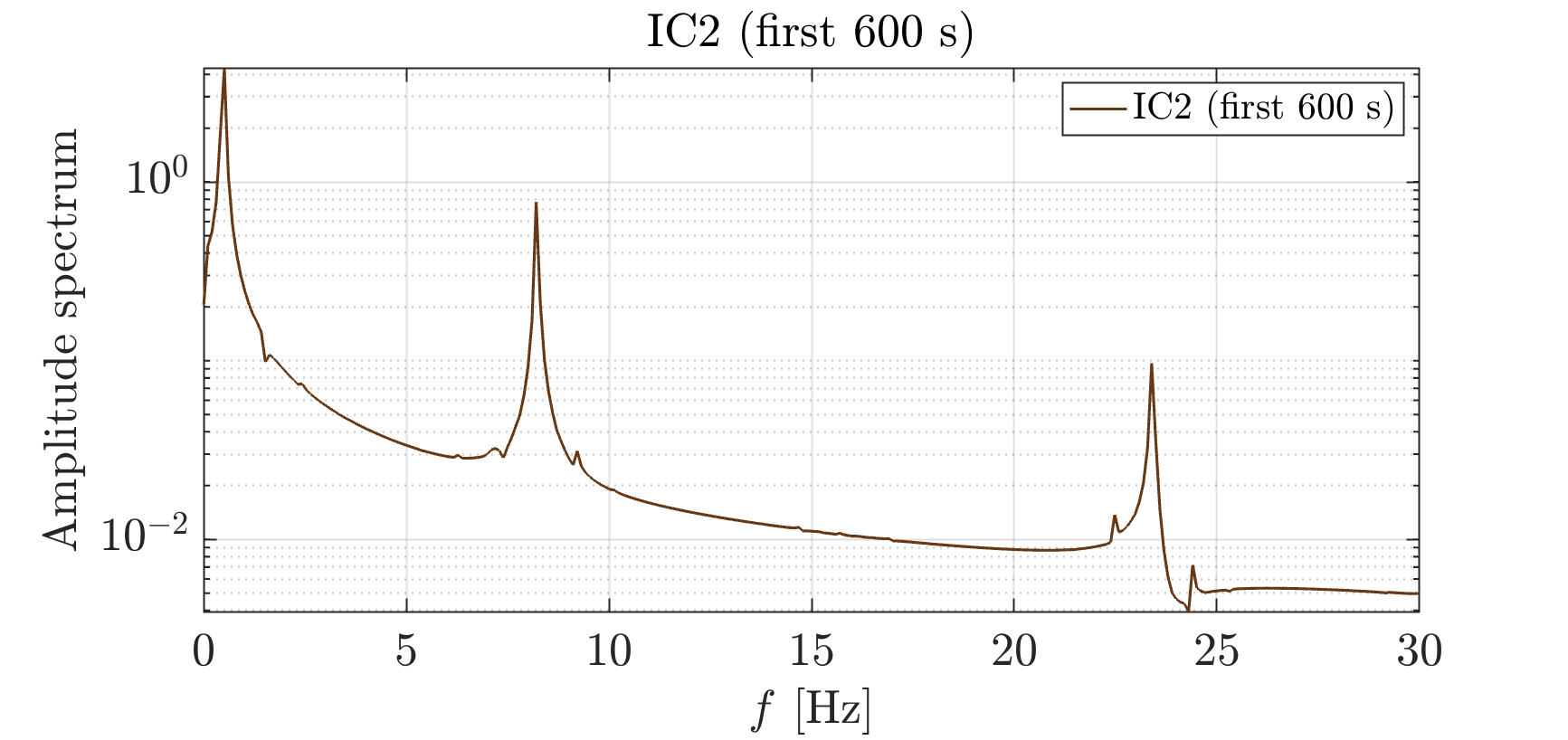}
  \caption{Amplitude spectrum obtained from the fast Fourier transform of the initial transient of full model trajectory of the Furuta pendulum initialized at IC2, confirming three distinct peaks associated with the unstable 3-torus predicted by the parametric SSM-reduced model.}
  \label{fig:Trajjj}
\end{figure}

The unstable 3-torus predicted by the parametric SSM-reduced model is verified in the full system in Fig.~\ref{fig:Trajjj}, which shows the amplitude spectrum of the initial transient of the trajectory initialized at the brown initial condition. The spectrum reveals three distinct peaks, corresponding to the 3-torus near which the initial condition lies.  

\subsection{Time delay and spatial discretization: microchaotic attractor modeling from numerical and experimental data}

We now discuss the data-driven modeling of an experimentally observed chaotic attractor arising from the spatial discretization associated with stabilizing the origin under delayed PD control. We start by extending the previous analysis by incorporating quantization modeling to represent spatial discretization. The controller parameters are unchanged from those used before as in Fig.~\ref{fig:stabMap}, namely $(K_{\rm P}, K_{\rm D})=(15.5\,\,\mathrm{V},5.45\,\,\mathrm{Vs})$ with $\Delta t=25\,\,\mathrm{ms}$ for the numerical simulations, while discretization is modeled via a floor operator.
For experiments, the parameters are $(K_{\rm P}, K_{\rm D}) = (25\,\mathrm{V},\,6\,\mathrm{Vs})$ and $\Delta t = 20\,\mathrm{ms}$, with a realistic quantization step matching the hardware resolution.

The first step is to estimate the fractal dimension of the chaotic attractor from experimental data. 
The estimated value is $4.39\pm 0.06$. For an SSM to contain such an attractor, its dimension must exceed this value. 
Given the oscillatory nature of the system, the eigenvalues of the linearized system at the equilibrium, occur in complex conjugate pairs. Therefore, the dimension of a mixed-mode SSM containing the attractor must be an even number larger than $4.39\pm 0.06$. We select $d=6$, i.e., construct a 6D mixed-mode SSM from experimental data.  

We project experimental trajectories onto the reduced linear subspace in the delay-embedding space using the optimization-based procedure described in \cite{Cenedese2022}, yielding reduced coordinates for training. In this reduced space, the discrete mapping generating the SSM-reduced dynamics is approximated via radial basis functions (RBFs). The RBFs are then used to forecast the state at the next time step, while the SSM parametrization lifts the reduced state back to the physical coordinates.  

The two experimental trajectories of longest duration are used to train the SSM-based model, and one experimental test trajectory is reserved for testing. Each trajectory consists of a time series of the scalar observable $\theta$. The use of only two training trajectories is sufficient, provided that they are of adequate duration, as we leverage the topological transitivity of the chaotic attractor. As a result, the two trajectories jointly explore the entire invariant set in a representative and informative manner.

A directly calculated prediction error along the test trajectory is not a good measure of model accuracy in the presence of a chaotic attractor, due to the sensitive dependence of the dynamics on initial conditions: any initial disturbance, however small, will be amplified by the system’s flow map due to a positive Lyapunov exponent. This leads to inevitable divergence between predicted and observed trajectories over time, even under a highly accurate approximation of the ODE governing the SSM-reduced dynamics. Therefore, a reduced model should not be expected to track individual chaotic trajectories indefinitely, but rather to reproduce long-term statistical properties of the dynamics. Accordingly, we assess the performance of the SSM-reduced model by comparing its ability to reproduce the leading Lyapunov exponent and the probability density distributions of the chaotic attractor.  

\subsubsection{Numerical simulations}

Following the modeling of quantization error proposed in \cite{Stepan2017, Haller1996}, the expected microchaotic attractor emerges in the vicinity of the fixed point to be stabilized. As illustrated in Figure~\ref{fig:NumSimPolyCombined}, 3D projections suggest that the attractor remains close to the stable 2-torus of the delayed, but spatially not discretized, counterpart system. Numerical simulations and experiments indicate that the quantization error transforms the stable 2-torus into a saddle-type 2-torus. Specifically, initial conditions placed near the original 2-torus are transiently attracted towards its remnant before being repelled along its unstable manifold, ultimately converging to a microchaotic attractor. During this transient regime, whose duration depends on the initial condition, the dynamics display quasiperiodic motion characterized by two interacting frequencies, after which the trajectories become chaotic.

At this stage, we first approximate the SSM-reduced dynamics using polynomials. Since trajectory prediction is unreliable for chaotic systems due to their sensitivity to initial conditions, our evaluation focuses on long-term statistics. Probability density functions (PDFs) of the reduced coordinates, shown in Figure~\ref{fig:NumSimPolyCombined}, demonstrate that these polynomial-based SSM models generate a smooth (in fact, analytic) 6D ODE system that has a smooth chaotic attractor. This attractor provides a close and smooth approximation of the non-smooth chaotic attractor of the full system. However, this polynomial approximation does not capture the leading Lyapunov exponent with sufficient accuracy: it yields the value $0.079$ in contrast to $0.048$ obtained from the full system's numerical simulations. This discrepancy motivates us to use radial basis functions with linear kernel in the subsequent section, which resulted in improved statistical predictions in previous studies (see, e.g., \cite{Liu2024, Kaszs2024, Xu2024, KaszasHaller2025}).

\newcommand{\imgEwidth}{0.55\textwidth} 
\newcommand{\imgFwidth}{0.45\textwidth} 

\begin{figure}[H] 
  \centering
  \begin{minipage}[c]{\imgEwidth}
    \centering
    \includegraphics[width=\linewidth]{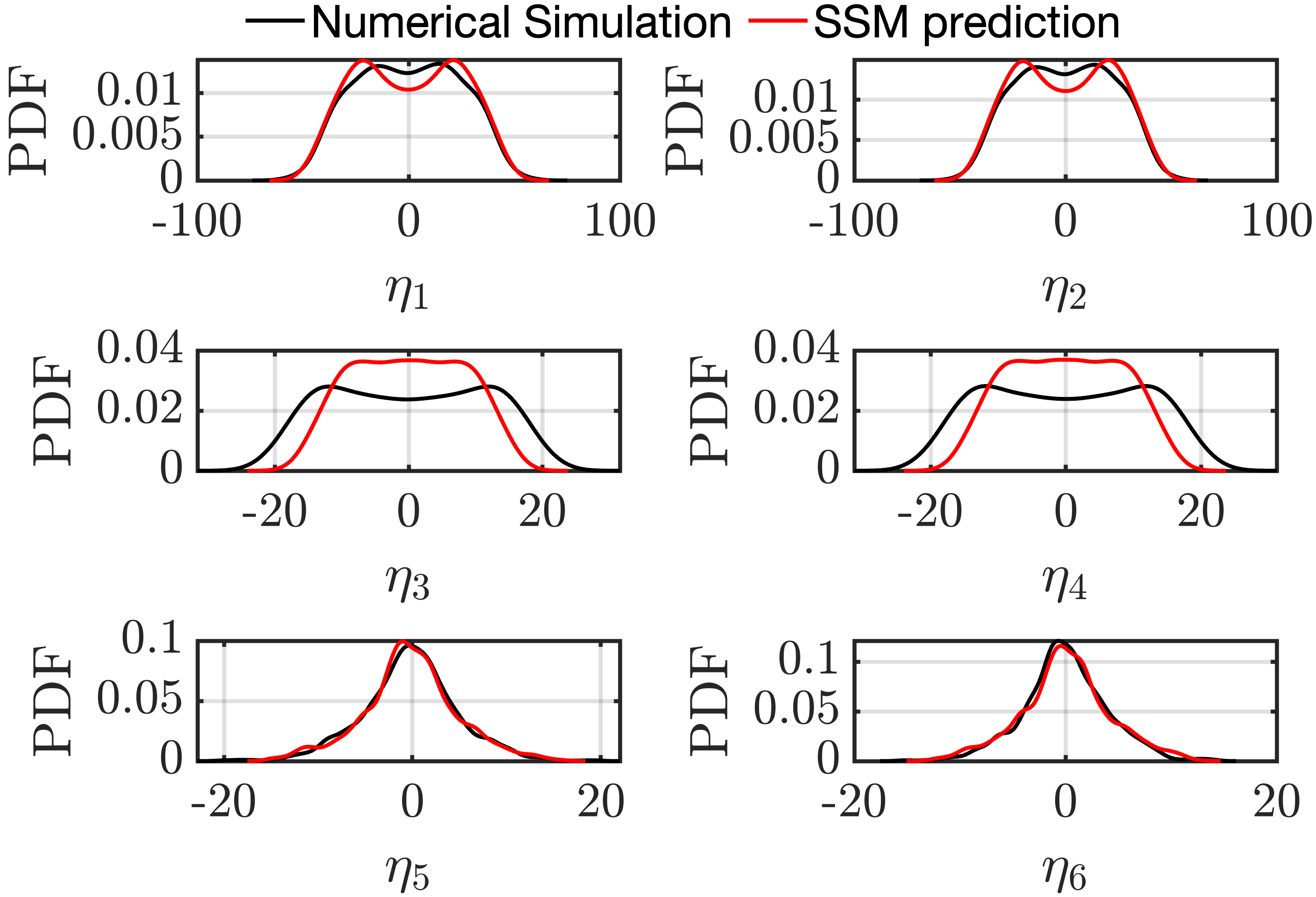}
  \end{minipage}\hfill
  \begin{minipage}[c]{\imgFwidth}
    \centering
    \includegraphics[width=\linewidth]{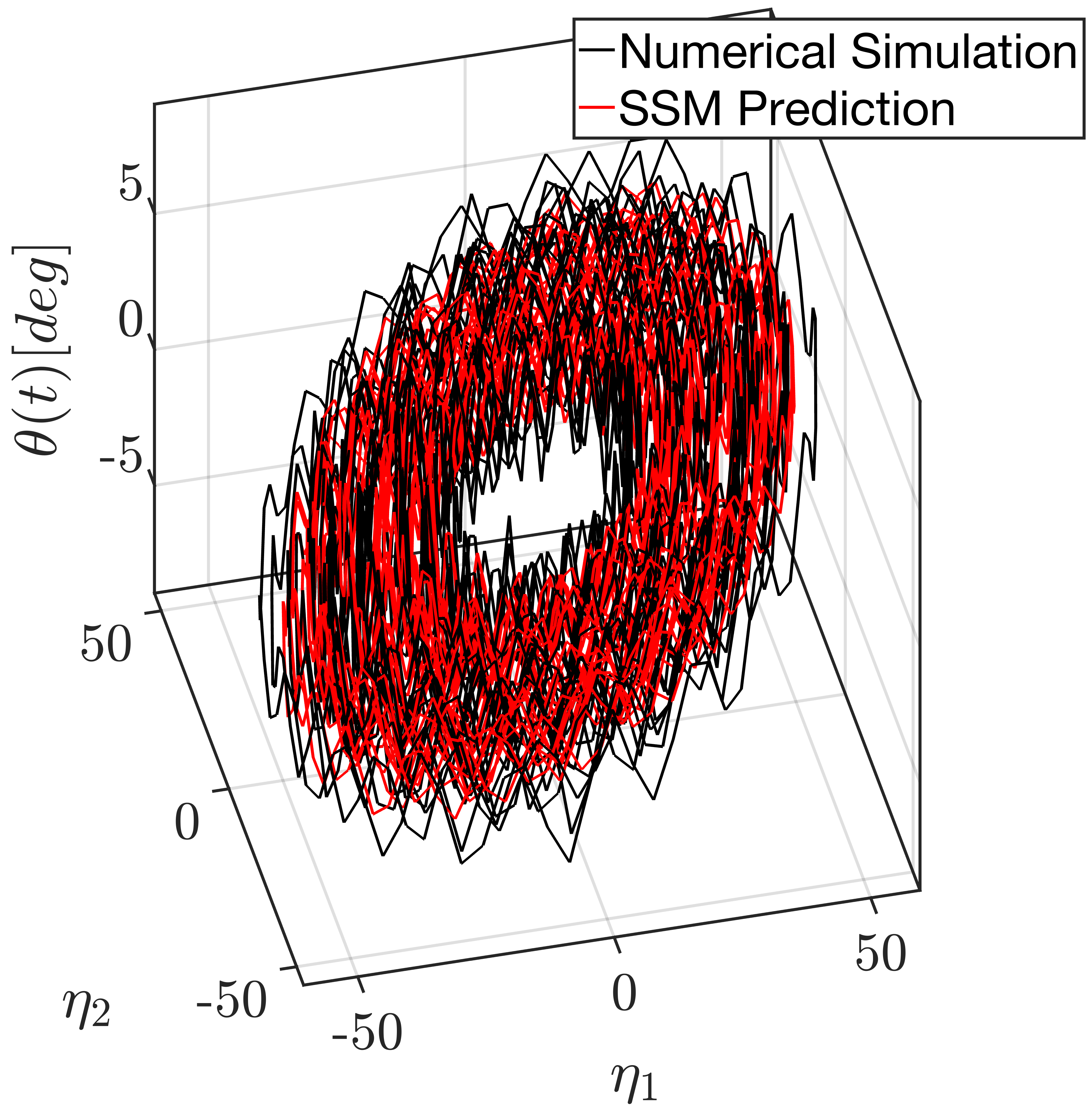 }
  \end{minipage}
  \caption{Left: Comparison of probability density functions of reduced coordinates obtained from numerical simulations (black) and the polynomial-based SSM model (red). Right: Chaotic attractor reconstructed from numerical simulations (black) and from the polynomial-based SSM model (red).}
  \label{fig:NumSimPolyCombined}
\end{figure}

\subsubsection{Experiments}

Next we use RBFs with linear kernel to approximate the discrete mapping governing the SSM-reduced dynamics in the available experimental data. The resulting model yields accurate results: the largest Lyapunov exponent is reproduced with a relative error of approximately $1.5\%$ ($0.132$ from the model versus $0.134$ from the experiments). In addition, the probability density functions of the reduced coordinates, shown in Figure~\ref{fig:ExpRBFCombined}, indicate close agreement between experiments and SSM-model predictions.

\newcommand{\imgGwidth}{0.55\textwidth} 
\newcommand{\imgHwidth}{0.45\textwidth} 

\begin{figure}[H] 
  \centering
  \begin{minipage}[c]{\imgGwidth}
    \centering
    \includegraphics[width=\linewidth]{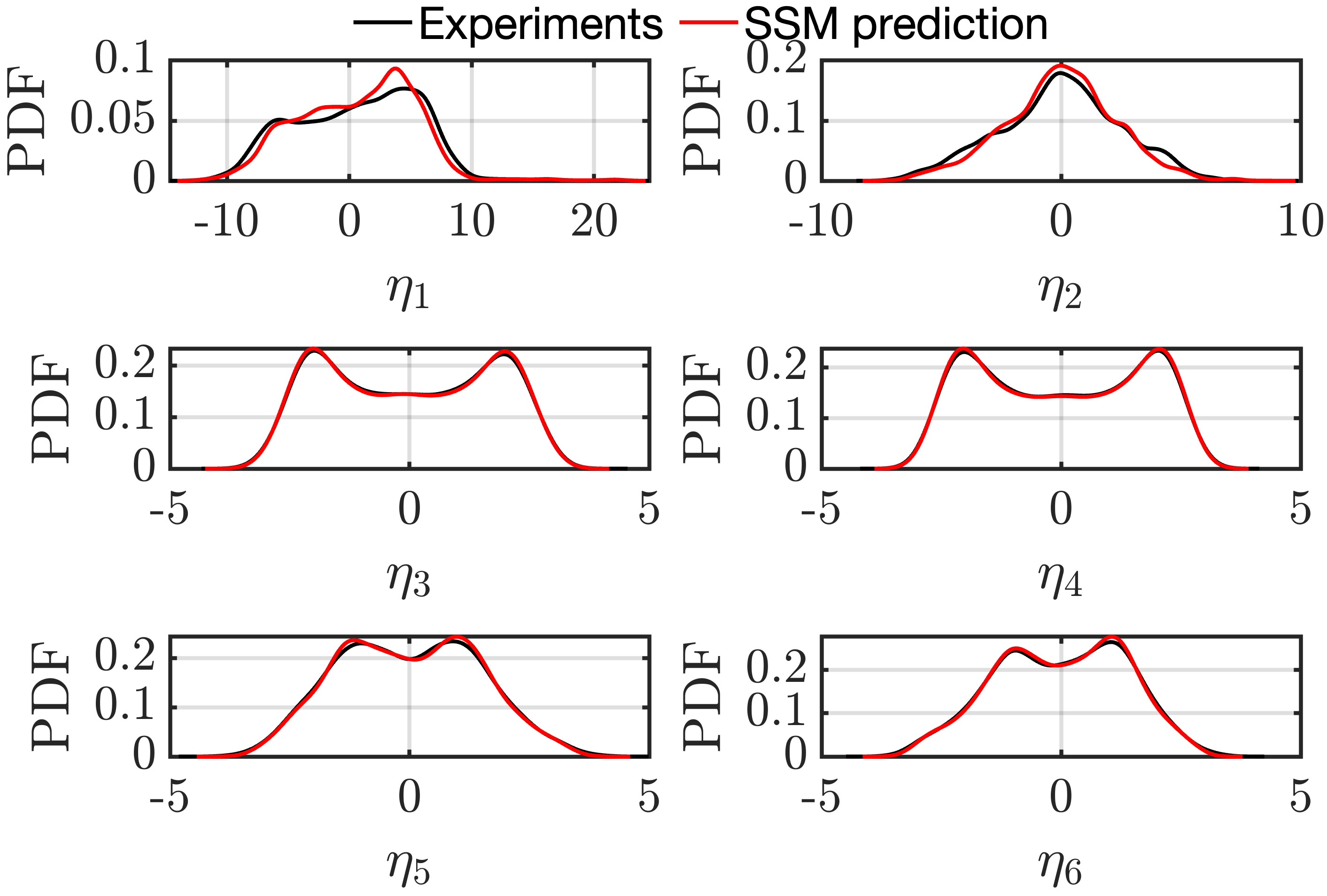}
  \end{minipage}\hfill
  \begin{minipage}[c]{\imgHwidth}
    \centering
    \includegraphics[width=\linewidth]{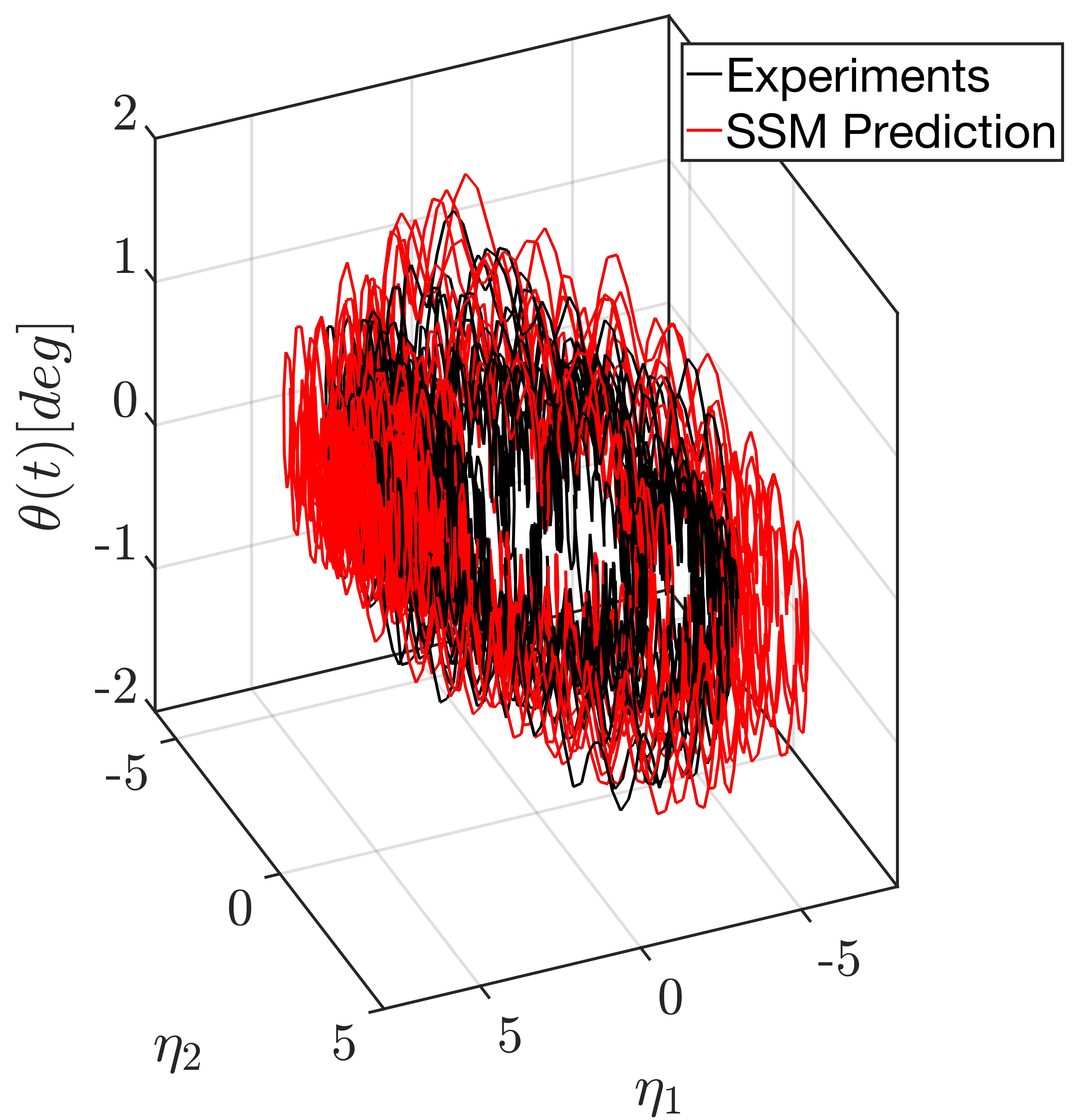}
  \end{minipage}
  \caption{Left: Comparison of probability density functions of reduced coordinates from experimental data (black) and the RBF-based SSM model (red). Right: Chaotic attractor reconstructed from experimental data (black) and from the RBF-based SSM model (red).}
  \label{fig:ExpRBFCombined}
\end{figure}

Trajectory prediction and the corresponding frequency spectrum are shown in Figure ~\ref{fig:TrajectoryPred123}, comparing experimental data with SSM-based model predictions. The spectrum exhibits two dominant peaks, consistent with the observation that the initial transient remains close to the remnants of a two-dimensional torus. In contrast, the trajectory comparison illustrates the strong sensitivity to initial conditions: even minimal discrepancies, likely introduced during the projection and subsequent lifting of the experimental initial condition to initialize the SSM model, result in a rapid divergence between the two trajectories.

\newcommand{\imgIwidth}{0.48\textwidth} 
\newcommand{\imgLwidth}{0.48\textwidth} 

\begin{figure}[H]
  \centering
  \begin{minipage}[c]{\imgIwidth}
    \centering
    \includegraphics[width=\linewidth]{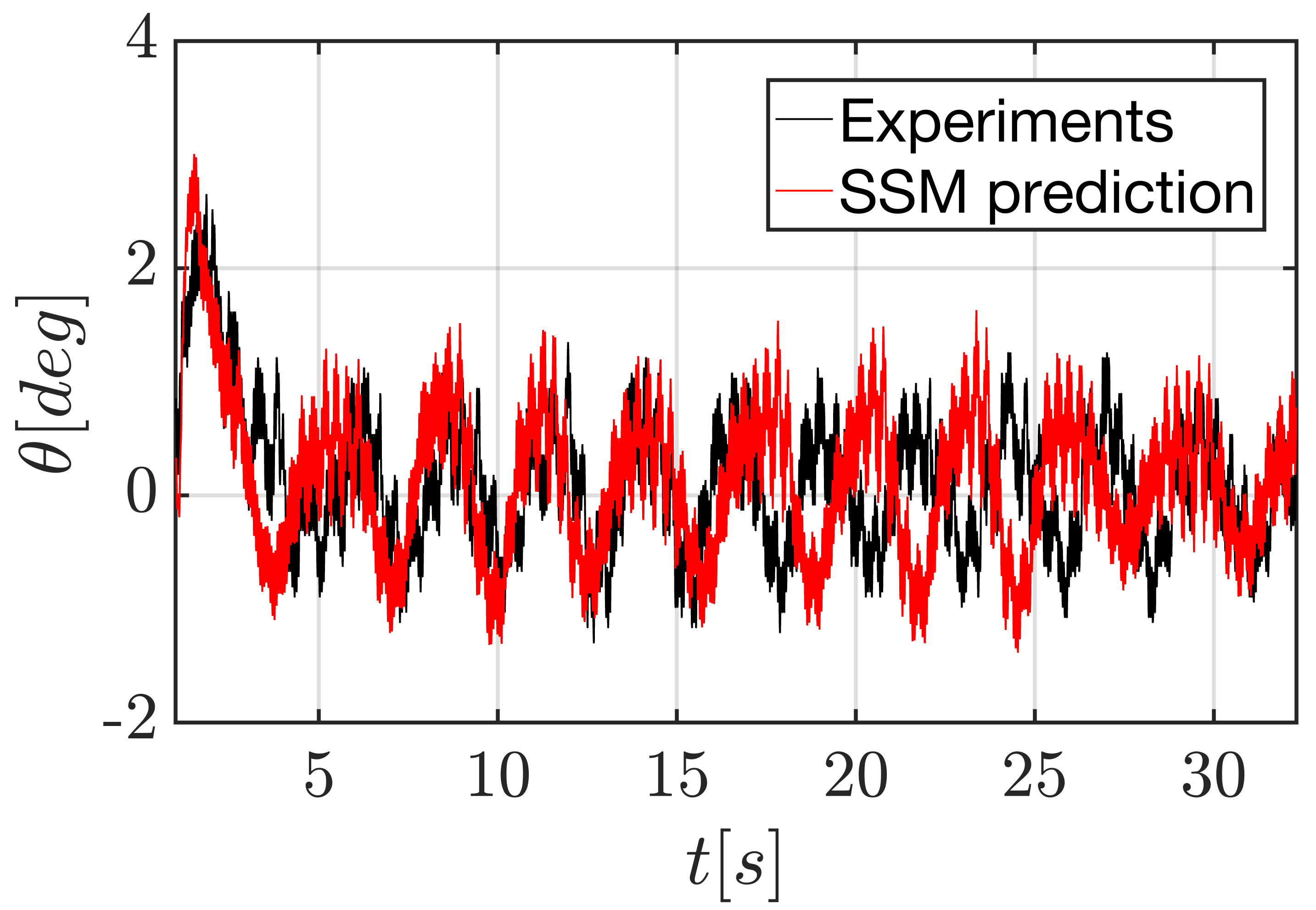}
  \end{minipage}\hfill
  \begin{minipage}[c]{\imgLwidth}
    \centering
    \includegraphics[width=\linewidth]{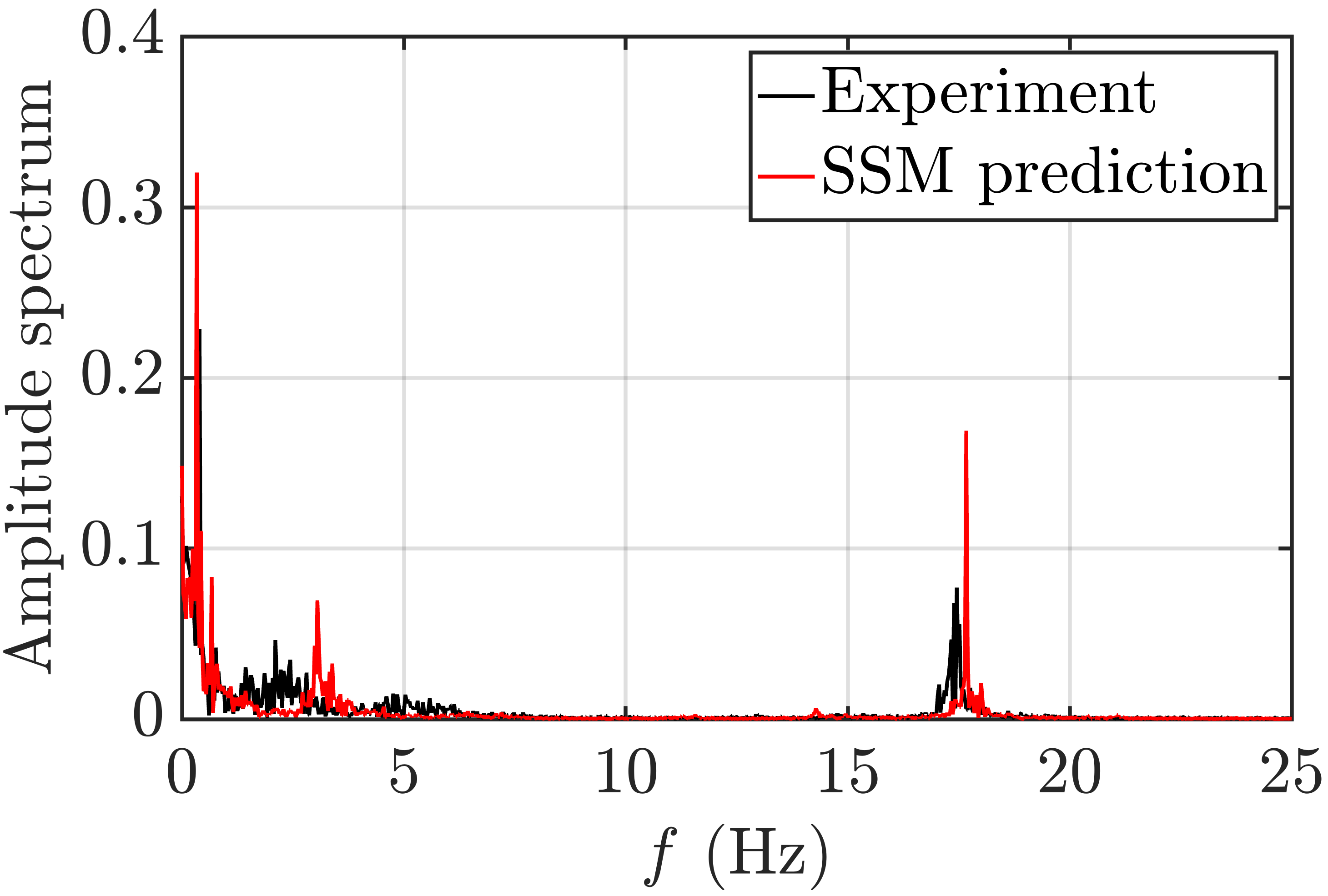}
  \end{minipage}
  \caption{Left: Time series of the scalar observable $\theta$ for a test trajectory. Right: Amplitude spectrum obtained via FFT of the same trajectory: experimental data (black) compared with SSM-based predictions (red).}\label{fig:TrajectoryPred123}

\end{figure} 

\section{Conclusions}

We have shown that a data-driven, SSM-based reduction captures the nonlinear and delay-induced dynamics of the Furuta pendulum subject to time-delayed PD controllers with spatial discretization. Specifically, exploiting the smooth dynamics of the uncontrolled system and the piecewise-analytic control action induced by a Zero-Order Hold logic, we constructed and validated a family of autonomous SSMs parametrized by the feedback delay.  

Numerical simulations of the Furuta pendulum showed that a 4D, data-driven, parametric SSM model accurately captures the nearby steady states, as well as local and global bifurcations of the full mechanical system. In particular, the model reproduces the transition from a stable two-torus, which explains the observed quasi-periodic motion of the pendulum, together with two saddle-type limit cycles. This transition occurs through a heteroclinic bifurcation that gives rise to an unstable three-torus and introduces a new route to instability, leading to the pendulum falling down. The model also captures a subsequent subcritical Hopf bifurcation. In this bifurcation, the three-torus shrinks and the two-torus loses stability, leaving only the $\mathcal{O}$(1)-distant fixed points introduced by the controller as attractors corresponding to tilted pendulum sustained rotation. The SSM model retained predictive accuracy at parameter values not used in the training, with phase errors evolving as expected from the intrinsic sensitivity of dynamical systems to perturbations in initial conditions. To support model validation, we introduced an error metric based on Dynamic Time Warping, which distinguishes between phase and amplitude errors and provides a more informative measure of model validation.  

We have also carried out SSM-reduced modeling on experimental data of the Furuta pendulum that incorporated spatial discretization and emergent microchaos. The fractal dimension of the experimental chaotic attractor was approximately 4.4, motivating the construction of a 6D SSM. Reduced dynamics in the form of a discrete map were identified from two experimental training trajectories using radial basis functions with linear kernel. The resulting model successfully reproduced short-term test trajectories and accurately captured key long-term statistical properties, including the leading Lyapunov exponent and the probability density functions of the reduced coordinates.  

Overall, these results establish SSM reduction as a viable and robust tool for system identification, model reduction, as well as local and global bifurcation prediction, under delay and nonsmooth effects. The approach described here applies more broadly to PID-stabilized mechanical systems with time-delayed feedback and spatial discretization, which can be reformulated into autonomous systems as illustrated here.
\section{Aknowledgments}
We are grateful to Bálint Kaszás, Roshan Kaundinya, Leonardo Bettini for useful comments on this work, and to James King for useful comments on the code. This work was supported by the Swiss National Science Foundation (Grant No. $200021\_214908$) and by the Hungarian National Research, Development and Innovation Office (NKFI-KKP-133846 and EKÖP-25-3-BME-138).

\section{Data availability}
The software and data supporting the findings of this study are available at
\url{https://github.com/haller-group/SSMLearn/tree/main/examples/parametricSSM_TimeDelay}.
\newpage
\appendix
\section{Data-driven spectral submanifolds for autonomous systems}\label{SectionAppA}

The theoretical foundations of Spectral Submanifolds (SSMs) in autonomous nonlinear dynamical systems, as well as their reconstruction via delay embedding from data, have been extensively developed in prior works (see \cite{Haller2025} for an introduction and overview). We briefly recall the most relevant results to the present analysis in this section.

\subsection{SSMs in autonomous nonlinear dynamical systems}
Consider an autonomous nonlinear system of the form
\begin{equation}\label{eq:1}
\dot x = A x + f_0(x),\quad x\in\mathbb R^n,\quad A\in\mathbb R^{n\times n},\quad f_0(x)=\mathcal O\bigl(|x|^2\bigr),\quad f_0\in C^r,\quad r\ge1,
\end{equation}
where $r \in \mathbb{N}^+ \;\cup\;\{\infty,a\}$. Assume that the eigenvalues of $A$ follows:
\begin{equation}\label{eq:stabAss}
Re \lambda_j < 0,\quad j = 1,\dots,n,
\end{equation}
which means that the origin is an asymptotically stable fixed point. Similar findings can be obtained for unstable fixed points (with $Re \lambda_j > 0$, $j = 1,\dots,n$), by inverting the direction of time.
A spectral subspace is the direct sum of arbitrary eigenspaces $E_{j_1,\dots,j_q}$ of $A$, and its smoothest nonlinear continuation is what we call the (primary) Spectral Submanifold for the nonlinear system~\eqref{eq:1}. 
\\
For a $q$‑dimensional spectral subspace $E$, we define the following linear combination of eigenvalues:
\[
\langle m,\lambda\rangle_E := m_1\lambda_{j_1} + \cdots + m_q\lambda_{j_q},
\quad
\lambda_{j_k}\in \mathrm{Spect}(A|_E),
\quad
m\in\mathbb{N}^q,
\]
with the order of the non-negative integer vector $m$ defined as $|m| := m_1 + \cdots + m_q$.
The relative spectral quotient \(\sigma(E)\) is defined for any spectral subspace \(E\) of the linear operator \(A\) as
\[
\sigma(E)
\; \doteq \;
\mathrm{Int}\!\Biggl[\,
\frac{\displaystyle\min_{\lambda\in \mathrm{Spect}(A)\setminus \mathrm{Spect}(A\lvert_E)} Re\lambda}
     {\displaystyle\max_{\lambda\in \mathrm{Spect}(A\lvert_E)} Re\lambda}
\Biggr].
\]
\noindent
The primary findings of significance in this study from \cite{Haller2025, Haller2016} can be summarized as follows: \\ \textbf{Theorem 1. } Consider a spectral subspace $E$ under the assumption that the low-order nonresonance conditions
\begin{equation}\label{eq:nonResCond}
\langle m,\lambda\rangle_E \neq \lambda_\ell,
\quad
\lambda_\ell \notin \mathrm{Spect}(A|_E),
\quad
2 \le |m| \le \sigma(E)
\end{equation}
hold for all eigenvalues $\lambda_\ell$ of $A$ that lie outside the spectrum of $A|_E$. Then the following statements are valid:
\begin{enumerate}[label=(\roman*)] 
    \item A class $C^r$ SSM, referred to as the primary SSM $\mathcal{W}(E)$, exists and it is tangent to the spectral subspace $E$ at $x = 0$. Additionally, it follows that $\dim \mathcal{W}(E) = \dim E$. 
    \item The manifold $\mathcal{W}(E)$ is unique among all $C^{\sigma(E)+1}$ invariant manifolds with the properties outlined in (i).
  \item If $f_0$ is jointly $C^r$ in $x$ and an additional parameter vector $\mu$, then the SSM $\mathcal{W}(E)$ is jointly $C^r$ in $x$ and $\mu$. Specifically, if $f_0(x,\mu)$ is either $C^\infty$ or analytic, the manifold $\mathcal{W}(E)$ persists under small perturbations in the parameter $\mu$ and will depend on these perturbations in a $C^\infty$ or analytic fashion, respectively.
  \item The primary SSM $\mathcal{W}(E)$ and the dynamics within $\mathcal{W}(E)$ depend smoothly on both $x$ and $\mu$. 
\end{enumerate}

\bigskip
Theorem 1 is also applicable when the $x=0$ fixed point is hyperbolic but not necessarily stable and the spectrum of $A$ contains no resonance. This means that
\begin{equation}
\mathrm{Re}\,\lambda_j \neq 0,\quad j=1,\dots,n, \quad \langle m,\lambda\rangle \neq 0, \quad \lambda \in \mathrm{Spect}(A), \quad  |m| >2,
\end{equation}
hold instead of~\eqref{eq:stabAss}-\eqref{eq:nonResCond} and $r=\infty$ is satisfied (see \cite{Haller2025} for more details).

\bigskip

Theorem 1 establishes a structured hierarchy $\mathcal{W}(E_{1}) \;\subset\; \mathcal{W}(E_{2}) \;\subset\;\dots\;\subset\; \mathcal{W}(E_{n-1})$ of SSMs which serves as an invariant nonlinear continuation of the spectral subspace hierarchy \(E_{1}\subset E_{2}\subset \dots\subset E_{n-1}\). In applications, the dominant dynamics is captured by the slowest member, \(\mathcal{W}(E_{1})\), of the SSMs family. To obtain predictions for shorter time scales, one can systematically enlarge the SSM-reduced model by including \(\mathcal{W}(E_{k})\) with increasing \(k\).

\subsection{SSM reconstruction via delay embedding}
It is uncommon in practice to find physical systems whose state variables can all be observed simultaneously \cite{Cenedese2022, Liu2024}. This requires the application of invariant-manifold reconstruction techniques that depend solely on a restricted set of observables. Utilizing the Takens' delay‐embedding theorem \cite{Takens1981}, one can effectively embed an invariant manifold of dimension $d$ within the space of $2d+1$ (or greater) time‐shifted samples derived from a generic scalar observable, as recently explored in \cite{Axs2023}. The reduced dynamics 
\[
\dot{\eta} = R(\eta), \qquad \eta \in \mathbb{R}^d,
\]
defined on a $d$-dimensional SSM, $\mathcal{W}(E)$, with flow map 
$R^t \colon \eta_0 \mapsto \eta(t)$, 
is conjugate to the restriction of the full flow map 
$F^t \colon \mathbb{R}^n \to \mathbb{R}^n$ 
of system~\eqref{eq:1} to 
$\mathcal{W}(E) \subset \mathbb{R}^n$, 
for a smooth parametrization 
$M \colon \mathbb{R}^d \to \mathbb{R}^n$ 
of the manifold $\mathcal{W}(E)$, i.e.,
\begin{equation}
    F^t \circ M = M \circ R^t.
\end{equation}

For a time-series $s(t) = \xi(x(t))$ obtained from observing system~\eqref{eq:1}, with $\xi \colon \mathbb{R}^n \to \mathbb{R}$, a delay-coordinate map with $m$ delays,
$\Psi \colon \mathbb{R}^n \to \mathbb{R}^m$,
is defined by aggregating $m$ uniformly sampled entries of the time series $s(t)$, with sampling interval $\Delta t$, into a vector,
\begin{equation}
    y(t) = \Psi(x(t))
    = [\,s(t),\, s(t+T_s),\, \dots,\, s(t+(m-1)T_s)\,]^\mathrm{T}
    \in \mathbb{R}^m.
\end{equation}

The map $\Psi$ is a function that transforms points in the phase space of system ~\eqref{eq:1} into points in the delay-embedded space $\mathbb{R}^m$.  Let $F_\Psi\colon\mathbb{R}^m\to\mathbb{R}^m$ represent the flow map generated by the flow of system ~\eqref{eq:1} inside this delay-embedding space.  By definition, if $\xi(0)=0$, the equilibrium point at $x=0$ of system ~\eqref{eq:1} corresponds to a fixed point at $y=0$, i.e. $F_\Psi(0)=0$. By the Takens embedding theorem \cite{Takens1981} for a generic observable function $\xi$, under nondegeneracy conditions on the SSM-reduced dynamics, if $m>2d$, then the delay-coordinate map $\Psi$ confined to $\mathcal{W}(E)$ embeds $\mathcal{W}(E)$ into $\mathbb{R}^m$ with probability one. Consequently, $\widetilde{\mathcal{W}}=\Psi(\mathcal{W}(E))$ is diffeomorphic to the manifold $\mathcal{W}(E)$.

The theoretical foundation for identifying an SSM attached to the $y=0$ origin of the $m$-dimensional delay-embedding space is provided by the conjugation of the delay-embedded dynamics $F_\Psi$ on $\widetilde{\mathcal{W}(E)}$ and the reduced dynamics restricted to the invariant manifold $\mathcal{W}(E)$
\begin{equation}
F_\Psi\circ\Psi = \Psi\circ F^t\bigl\lvert_{\mathcal{W}(E)}.
\end{equation}

\subsection{Data‐driven modeling of the spectral submanifold in a delay embedding space}
Following \cite{Liu2024}, let the coordinate $y\in\mathbb{R}^m$ represent data acquired from system ~\eqref{eq:1} by direct observations and subsequent delay‐embedding.  Let $\eta\in\mathbb{R}^d$ represent the reduced coordinates corresponding to the spectral subspace $E$ of $D\Psi(E)$.  Locally, the embedded SSM can be approximated by a polynomial expansion in the vicinity of the anchor fixed point, as a graph tangent to this spectral subspace. 

Consequently, calling $\mathcal{K}$ the maximum polynomial order of the expansion, we seek a parametrization of $\mathcal{W}(E)$ as

\begin{equation}\label{eq:M_eta}
M(\eta)
= V_1\,\eta + \sum_{|\mathbf{k}|=2}^{\mathcal{K}} V_{\mathbf{k}}\,\eta^{\mathbf{k}},
\qquad 
\mathbf{k}=(k_1,\dots,k_d)\in\mathbb{N}^d,\quad
V_1,V_{\mathbf{k}}\in\mathbb{R}^{m\times d},
\end{equation}
where $\eta^{\mathbf{k}}:=\eta_1^{k_1}\eta_2^{k_2}\cdots\eta_d^{k_d}$ denotes the multivariate monomial corresponding to the multi‐index $\mathbf{k}$.

Reduced coordinates on $\mathcal{W}(E)$ can be defined as $\eta = V_1^\mathrm{T} y$, where the matrix $V_1\in\mathbb{R}^{m\times d}$ has orthonormal columns that include the tangent space of the manifold $\widetilde{\mathcal{W}(E)}(E)\subset\mathbb{R}^m$ at the point $y=0$. We use SSMLearn \cite{Cenedese2022} to solve the optimization problem

\begin{equation}\label{eq:Vopt}
V^* = [\,V_1^*,\,V_{2:\mathcal{K}}^*\,]= \argmin_{V_1,\,\{V_{\mathbf{k}}\}_{|\mathbf{k|}\ge2}}
   \sum_j
   \Big\lVert
     y_j 
     - V_1\,V_1^\mathrm{T} y_j 
     - \sum_{|\mathbf{k}|=2}^{\mathcal{K}} V_{\mathbf{k}}\,(\eta_j)^{\mathbf{k}}
   \Big\rVert_2^2  \qquad 
\end{equation}
\newline
subject to $V_1^\mathrm{T} V_1 = I, \quad V_1^\mathrm{T} V_{\mathbf{k}} = 0$ for all $\mathbf{k}$ with $|\mathbf{k}|\ge2$.
Here, the second constraint represents a fundamental nonlinear extension of the principal component analysis.

\subsection{Data‐driven modeling of the SSM-reduced dynamics}
One can employ two alternative ways of approximating the reduced dynamics within the SSM: polynomial regression and radial basis function interpolation (RBF).

\subsubsection*{Polynomial reduced dynamics}
We approximate the reduced dynamics as
\begin{equation}\label{eq:reduced-dynamics}
\dot{\eta} 
= R_{1}\,\eta 
+ \sum_{|\mathbf{k}|=2}^{\mathcal{K}} R_{\mathbf{k}}\,\eta^{\mathbf{k}},
\qquad 
\mathbf{k}=(k_1,\dots,k_d)\in\mathbb{N}^d,
\quad
R_{\mathbf{k}}\in\mathbb{R}^{d\times d},
\end{equation}
where the coefficient matrices are obtained via regression as

\begin{equation}\label{eq:regression}
R^{*}
=\argmin_{R}
\sum_{j}
\Bigl\|
\dot{\eta}_{j} 
-\sum_{|\mathbf{k}|=1}^{\mathcal{K}} R_{\mathbf{k}}\,(\eta_{j})^{\mathbf{k}}
\Bigr\|_2^2,
\qquad 
\mathbf{k}=(k_1,\dots,k_d)\in\mathbb{N}^d, \quad R = \bigl[\,R_{1},\,R_{2:\mathcal{K}}\,\bigr]
\;\in\;
\mathbb{R}^{d\times\sum_{k=1}^{\mathcal{K}}d_{k}}.
\end{equation}
The number of distinct multivariate monomials 
of total degree~$k$ is given by $d_k = \binom{d+k-1}{k}$. Subsequently, we diagonalize the linear part of the reduced dynamics~\eqref{eq:reduced-dynamics} 
through a linear change of coordinates, $\xi = W^{-1}\,\eta$, where $\xi \in \mathbb{R}^d$ denotes the modal coordinates and 
$W \in \mathbb{R}^{d\times d}$ is the matrix of eigenvectors of $R_{1}$.
In these new coordinates, the reduced dynamics~\eqref{eq:reduced-dynamics} take the form
\begin{equation}\label{eq:modal-ode}
\dot{\xi}
= \Lambda\,\xi
+ \sum_{|\mathbf{k}|=2}^{\mathcal{K}} N_{\mathbf{k}}\,\xi^{\mathbf{k}},
\qquad 
\mathbf{k} = (k_1, \dots, k_d) \in \mathbb{N}^d.
\end{equation}
Here, $\Lambda \in \mathbb{R}^{d\times d}$ is a diagonal matrix containing the eigenvalues of $R_{1}$, 
and each $N_{\mathbf{k}} \in \mathbb{R}^{d\times 1}$ collects the coefficients of the nonlinear terms 
associated with the multivariate monomial 
$\xi^{\mathbf{k}} := \xi_1^{k_1}\xi_2^{k_2}\cdots\xi_d^{k_d}$.

The coefficient matrices $N_{2:\mathcal{K}}$ will be defined by a recursive series of normal form transformations that maintain the near-resonant terms. This yields the extended normal form of the SSM-reduced dynamics as detailed in \cite{Cenedese2022, Haller2025}, expressed in polar coordinates. 

Classic normal form transformations of increasing order are typically only defined on progressively smaller neighborhoods of $y=0$, restricting the area of validity of our local SSM model. In contrast, extended normal forms remain valid on larger domains.

\newpage
\subsubsection*{RBF reduced dynamics}
We use radial basis functions with linear kernel, as in \cite{Xu2024}, to approximate the dynamics reduced to the SSM as a discrete mapping,
\begin{equation}
\eta_{n+1}
= F(\eta_n)
= \sum_{i=1}^K C_i\,k\bigl(\|\eta_n - \eta_i\|\bigr).
\end{equation}
Here \(k\) denotes a radial kernel function dependent upon the distance between \(\eta_n\) and the points \(\eta_i\) in the training dataset. We choose a linear kernel by defining \(k(r)=r\), for its performance and simplicity.

\subsection{Parametric SSM-reduced order models}
By the smooth dependence of the SSM and its reduced dynamics on the parameters, we can construct parameter-dependent SSMs. This is achieved by smoothly interpolating the coefficients of each monomial in the multivariate polynomial parametrization of both the manifold and its reduced dynamics (e.g. spline interpolation). Each coefficient is then obtained as a function of the parameter $\mu \in \mathbb{R}^p$. The fixed point at which the SSM is anchored must remain hyperbolic to ensure the smooth persistence of the SSM. This criterion enables the construction of reduced models that capture both local bifurcations of nearby steady states and global bifurcations.

As an example, the parameter-dependent extended normal form in a four-dimensional $(d=4)$ SSM constructed over oscillatory linear modes can be written in polar coordinates as

\begin{equation}
\begin{aligned}
\dot\rho_1\,\rho_1^{-1} &= a_1(\mu) + b_1(\mu)\,\rho_1^2 + c_1(\mu)\,\rho_2^2
  + \mathcal{O}\bigl(\lvert\rho\rvert^4\bigr),\\
\dot\rho_2\,\rho_2^{-1} &= a_2(\mu) + b_2(\mu)\,\rho_1^2 + c_2(\mu)\,\rho_2^2
  + \mathcal{O}\bigl(\lvert\rho\rvert^4\bigr),\\
\dot\theta_1 &= a_3(\mu) + \mathcal{O}\bigl(\lvert\rho\rvert^2\bigr),\\
\dot\theta_2 &= a_4(\mu) + \mathcal{O}\bigl(\lvert\rho\rvert^2\bigr).
\end{aligned}
\end{equation}

\subsection{Predictions for SSM containing chaotic attractors}

For a spectral submanifold (SSM) to contain a chaotic attractor, its dimension $d$ must exceed the fractal dimension of the attractor. In data-driven settings, wherein no prior knowledge of the attractor dimension is available, this fractal dimension is estimated from time series using the Grassberger–Procaccia (GP) algorithm \cite{GrassbergerProcaccia1983}. The GP method requires reconstructing trajectories on the attractor in a delay-embedding space, as described by \cite{Takens1981}.

In oscillatory systems, the SSM is further constrained to be even-dimensional, with $d$ taken as the smallest even integer greater than the estimated fractal dimension.  

Under spatial discretization~\eqref{eq:ODESD11}, we can construct a smooth SSM-reduced model from data to obtain an approximation of the non-smooth system. Specifically, we employ the open-source MATLAB package SSMLearn \cite{Cenedese2022} to compute the SSM parametrization, and reduced dynamics using radial basis functions (RBFs) with linear kernel, based solely on time series of a scalar observable.  

Our numerical experiments show that polynomials can provide a close smooth approximation of the non-smooth microchaotic attractor of the Furuta pendulum. However, in agreement with previous findings \cite{Xu2024, Liu2024}, we find that linear RBFs outperform polynomial SSMs approximations.
   
\section{PID problem setup}\label{SectionAppB}

Here we extend the derivation of Section~\ref{section3} to PID controllers by including the integral action.
\subsection{Uncontrolled system}

Consider the autonomous uncontrolled dynamical system,
\begin{equation}\label{eq:uncontrolledSysB}
    \dot x = f(x)=Ax + f_0(x), \quad x\in\mathbb{R}^n, \quad A\in\mathbb{R}^{n\times n}, \quad f_0 = \mathcal{O}(|x|^2)\in C^r, \quad r\in \mathbb{N}^+ \cup \{ \infty, a \} .
\end{equation}
where the origin $x=0$ is a hyperbolic fixed point to be stabilized by feedback.

\subsection{PID controlled system}
We define an observable vector
\[
\xi:\mathbb{R}^n \to \mathbb{R}^l,
\quad
\xi(0)=0,
\quad 
\xi \in C^1,
\]

\noindent
and the tracking error
\[
e(t) = \xi(0) - \xi(x(t)) = -\,\xi(x(t)),
\]
which vanishes at $x=0$. We employ a PID controller with gains $K_{\rm P},K_{\rm I},K_{\rm D}\in\mathbb{R}^{n\times l}$.
\noindent
Define the integral state
\[
x_{\rm I}(t)=\int_0^t \xi(x(s))\,\mathrm ds,
\qquad
\dot x_{\rm I}(t)=\xi(x(t)),
\]
so that the integral action is $u_{\rm I}(t)=-K_{\rm I}x_{\rm I}(t)$.

With the extended state
\[
\tilde x(t)=
\begin{pmatrix}
x(t)\\
x_{\rm I}(t)
\end{pmatrix}
\in\mathbb{R}^{n+l},
\]
the closed-loop dynamics read
\[
\dot{\tilde x}(t)=
\begin{pmatrix}
f(x(t))+u(x(t),x_{\rm I}(t))\\
\xi(x(t))
\end{pmatrix}.
\]
The closed-loop system is autonomous, with equilibrium $\tilde x=0$.

\subsection{PID controlled system with time-delayed feedback}\label{sec:subsec3.3}

Incorporating a delay $\tau>0$ in the feedback loop, the closed-loop dynamics satisfy the DDE
\begin{equation}\label{eq:delayPDcontrolledSys1B}
\dot{x}(t)
= f\bigl(x(t)\bigr)
- K_{\rm P}\,\xi\bigl(x(t-\tau)\bigr)
- K_{\rm I}\,x_{{\rm I}}(t-\tau)
- K_{\rm D}\,\frac{{\rm d}}{{\rm d}t}\,\xi\bigl(x(t-\tau)\bigr).
\end{equation}
In the extended $(n+l)$-dimensional state space
\begin{equation}
\dot{\tilde{x}}(t)
=
\begin{pmatrix}
\dot x(t) \\[4pt]
\dot x_{\rm I}(t)
\end{pmatrix}
=
\begin{pmatrix}
f\bigl(x(t)\bigr)
- K_{\rm P}\,\xi\bigl(x(t-\tau)\bigr)
- K_{\rm I}\,x_{\rm I}(t-\tau)
- K_{\rm D}\,\dfrac{{\rm d}}{{\rm d}t}\,\xi\bigl(x(t-\tau)\bigr)
\\[8pt]
\xi\bigl(x(t)\bigr)
\end{pmatrix}.
\label{eq:delayPIDB}
\end{equation}
\noindent We define the extended state vector as
\begin{equation}\label{eq:ExtensionB}
\tilde{x}_{\mathrm{ext}}(t)
=
\begin{pmatrix}
\tilde{x}_0(t) \\[2pt]
\tilde{x}_1(t) \\[2pt]
\tilde{x}_2(t) \\[2pt]
\vdots
\end{pmatrix},
\qquad
\tilde{x}_k(t) = 
\begin{pmatrix}
x_k(t) \\[2pt]
x_{{\rm I},k}(t)
\end{pmatrix}
=
\begin{pmatrix}
x(t-k\tau) \\[2pt]
x_{\rm I}(t-k\tau)
\end{pmatrix},
\quad
\tilde{x}_k \in \mathbb{R}^{\,n+l},
\quad
k \in \mathbb{N},
\end{equation}
The dynamics of the $k^{\rm th}$ layer evolve as
\begin{equation}
\dot{\tilde{x}}_k(t)
=
\begin{pmatrix}
\dot x_k(t) \\[4pt]
\dot x_{{\rm I},k}(t)
\end{pmatrix}
=
\underbrace{%
  \begin{pmatrix}
    f\bigl(x_k(t)\bigr) \\[6pt]
    \xi\bigl(x_k(t)\bigr)
  \end{pmatrix}}_{=\,F_k\bigl({x}_k(t)\bigr)}
\;+\;
\underbrace{%
  \begin{pmatrix}
    -K_{\rm P}\,\xi\bigl(x_{k+1}(t)\bigr)
    -K_{\rm I}\,x_{{\rm I},k+1}(t)
    -K_{\rm D}\,\dot{\xi}\bigl(x_{k+1}(t)\bigr) \\[6pt]
    0_{\,l\times1}
  \end{pmatrix}}_{=\,G_k\bigl(\tilde{x}_{k+1}(t)\bigr)}.
\end{equation}
Equivalently,
\begin{equation}
\dot{\tilde{x}}_k(t)
= F_k\bigl({x}_k(t)\bigr) + G_k\bigl(\tilde{x}_{k+1}(t)\bigr), \qquad  k=0,1,2,\dots \quad .
\end{equation} 
Note that each layer $k$ depends on the subsequent layer $k+1$. Finally, let us introduce the mappings
\[
F(\tilde{x}_{\mathrm{ext}}) \;=\;
\begin{pmatrix}
F_0({x}_0(t))\\[4pt]
F_1({x}_1(t))\\[4pt]
F_2({x}_2(t))\\[2pt]
\vdots
\end{pmatrix},
\qquad
G(\tilde{x}_{\mathrm{ext}}) \;=\;
\begin{pmatrix}
G_0(\tilde{x}_1(t))\\[4pt]
G_1(\tilde{x}_2(t))\\[4pt]
G_2(\tilde{x}_3(t))\\[2pt]
\vdots
\end{pmatrix}.
\]

Then
\begin{equation}\label{eq:DelayEquationB}
\dot{\tilde{x}}_{\mathrm{ext}}(t)
= F\bigl(\tilde{x}_{\mathrm{ext}}(t)\bigr)
+ G\bigl(\tilde{x}_{\mathrm{ext}}(t)\bigr).
\end{equation}
Equation~\eqref{eq:DelayEquationB} is an autonomous ODE representation of~\eqref{eq:delayPDcontrolledSys1B}, with equilibrium
$\tilde x_{\mathrm{ext}}=0$.

\subsection{PID controlled system with time-delayed feedback and spatial discretization}

We include spatial discretization due to quantized measurement and actuation. Define

\begin{equation}
x_{{\rm I},{\rm d}}(t) \;=\; \int_{0}^{t} h \left\lfloor \tfrac{1}{h}\,\xi\bigl(x(s)\bigr)\,\mathrm{d}s \right\rfloor,
\qquad
\dot{x}_{{\rm I},{\rm d}}(t) \;=\; h \left\lfloor \tfrac{1}{h}\,\xi(x(t)) \right\rfloor \text{ for almost every } t \in \mathbb{R},
\end{equation}

\noindent where $h>0$ is the quantization step. Including this discretization effect \cite{Stepan2017} yields

\begin{equation}\label{eq:realWorldSDB}
\begin{aligned}
\dot{\tilde{x}}_{\rm d}(t)
&=
\begin{pmatrix}
\dot{x}(t) \\[8pt]
\dot{x}_{{\rm I},{\rm d}}(t)
\end{pmatrix}=
\\[6pt]
&=
\begin{pmatrix}
f\bigl(x(t)\bigr) +
 h\left\lfloor  
   - K_{\rm P}\,\Bigl\lfloor \tfrac{1}{h}\,\xi \bigl(x(t-\rho(t))\bigr)\Bigr\rfloor
   - \tfrac{1}{h} K_{\rm I}\,x_{{\rm I},{\rm d}}\bigl(t-\rho(t)\bigr)
   - K_{\rm D}\,\Bigl\lfloor \tfrac{1}{h}\,\dot \xi\bigl(x(t-\rho(t))\bigr)\Bigr\rfloor
 \right\rfloor
\\[6pt]
h\Bigl\lfloor \tfrac{1}{h}\,\xi(x(t)) \Bigr\rfloor
\end{pmatrix},
\end{aligned}
\end{equation}

\noindent with $\tilde{x}_{{\rm d}}(t) = \bigl[x(t);\,x_{{\rm I},{\rm d}}(t)\bigr]$. 
Equation~\eqref{eq:realWorldSDB} models that the control action is computed and applied as an integer multiple of the quantization step $h$, and the time discretization follows the ZOH. 

Analogously to the construction in~\eqref{eq:ExtensionB}, we define the $k^{\rm th}$ layer as
\[
\tilde{x}_{{\rm d},k}(t) =
\begin{pmatrix}
x_{k}(t) \\[2pt]
x_{{\rm I},{\rm d},k}(t)
\end{pmatrix}
=
\begin{pmatrix}
x\!\bigl(t - \operatorname{sgn}(k)\,\rho(t)
           - H(k - 2)\,r (k-1)\,\Delta t \bigr) \\[4pt]
x_{{\rm I},{\rm d}}\!\bigl(t - \operatorname{sgn}(k)\,\rho(t)
           - H(k - 2)\,r (k-1)\,\Delta t \bigr)
\end{pmatrix}
\in \mathbb{R}^{\,n+l},
\quad k=0,1,2,\dots
\]
\noindent  We also introduce the extended state as
\begin{equation}\label{eq:Extension1B}
\tilde{x}_{\mathrm{ext},{\rm d}}(t)
=\begin{pmatrix}
\tilde{x}_{{\rm d},0}(t) \\[23pt]
\tilde{x}_{{\rm d},1}(t) \\[23pt]
\tilde{x}_{{\rm d},2}(t) \\[23pt]
\tilde{x}_{{\rm d},3}(t) \\[23pt]
\vdots
\end{pmatrix}
=
\begin{pmatrix}
x(t) \\[0pt]
x_{{\rm I},{\rm d}}(t) \\[10pt]
x(t-\rho (t)) \\[0pt]
x_{{\rm I},{\rm d}}(t-\rho (t)) \\[10pt]
x(t-\rho (t)-r\Delta t) \\[0pt]
x_{{\rm I},{\rm d}}(t-\rho (t)-r\Delta t) \\[10pt]
x(t-\rho (t)-2r\Delta t) \\[0pt]
x_{{\rm I},{\rm d}}(t-\rho (t)-2r\Delta t) \\[2pt]
\vdots
\end{pmatrix},
\end{equation}
At delay layer $k$, the dynamics takes the form
\begin{equation}\label{eq:delayed_layersB}
\begin{aligned}
\dot{\tilde{x}}_{{\rm d},k}(t)
&=
\begin{pmatrix}
\dot x_{k}(t) \\[2pt]
\dot x_{{\rm I},{\rm d},k}(t)
\end{pmatrix}=
\\[6pt]
&=
\underbrace{\begin{pmatrix}
f\bigl(x_{k}(t)\bigr) \\[6pt]
h \left\lfloor \tfrac{1}{h}\,\xi\bigl(x_{k}(t)\bigr) \right\rfloor
\end{pmatrix}}_{=\,F_{{\rm d},k}\bigl(\tilde{x}_{{\rm d},k}(t)\bigr)}
\;+\;
\underbrace{\begin{pmatrix}
h \Bigl\lfloor -K_{\rm P}\,\Bigl\lfloor \tfrac{1}{h}\,\xi\bigl(x_{k+1}(t)\bigr)\Bigr\rfloor
- \tfrac{1}{h}K_{\rm I}\,x_{{\rm I},{\rm d},k+1}(t)
- K_{\rm D}\,\Bigl\lfloor \tfrac{1}{h}\,\dot{\xi}\bigl(x_{k+1}(t)\bigr)\Bigr\rfloor \Bigr\rfloor \\[6pt]
0_{\,l\times 1}
\end{pmatrix}}_{=\,G_{{\rm d},k}\bigl(\tilde{x}_{{\rm d},k+1}(t)\bigr)\;=\;\text{input from next delay layer}}.
\end{aligned}
\end{equation}

\noindent As before, we extend the phase space and introduce the mappings
\[
F_{\rm d}(\tilde{x}_{\mathrm{ext},{\rm d}})
=
\begin{pmatrix}
F_{{\rm d},0}(x_{{\rm d},0}(t))\\[4pt]
F_{{\rm d},1}(x_{{\rm d},1}(t))\\[4pt]
F_{{\rm d},2}(x_{{\rm d},2}(t))\\[2pt]
\vdots
\end{pmatrix},
\qquad
G_{\rm d}(\tilde{x}_{\mathrm{ext},{\rm d}})
=
\begin{pmatrix}
G_{{\rm d},0}(\tilde{x}_{{\rm d},1}(t))\\[4pt]
G_{{\rm d},1}(\tilde{x}_{{\rm d},2}(t))\\[4pt]
G_{{\rm d},2}(\tilde{x}_{{\rm d},3}(t))\\[2pt]
\vdots
\end{pmatrix}.
\]
\noindent
We obtain the autonomous, infinite-dimensional ODE version of ~\eqref{eq:realWorldSDB}, incorporating temporal sampling, ZOH, feedback delay, and spatial discretization, with equilibrium at $\tilde{x}_{\mathrm{ext},{\rm d}}=0$:

\begin{equation}\label{eq:ODESD123B}
\dot{\tilde{x}}_{\mathrm{ext},{\rm d}}(t)
= F_{\rm d}\bigl(\tilde{x}_{\mathrm{ext},{\rm d}}(t)\bigr)
+ G_{\rm d}\bigl(\tilde{x}_{\mathrm{ext},{\rm d}}(t)\bigr).
\end{equation}

\let\clearpage\oldclearpage
\newpage
\addcontentsline{toc}{section}{References}
\printbibliography

\end{document}